\documentclass[a4,12pt]{amsart}
\usepackage{amsfonts,textcomp,amssymb,amsmath,amsthm}
\usepackage{color,ulem}
\usepackage{fullpage}
\usepackage{graphicx}
\usepackage[active]{srcltx}
\usepackage{mathabx}
\usepackage{dsfont,mathrsfs,stmaryrd,wasysym,amsbsy}
\usepackage{enumerate}

\newcommand{\EE}{\mathbb{E}}
\newcommand{\E}{\mathbb{E}}

\newcommand{\NN}{\mathbb{N}}
\newcommand{\nn}{\mathbb{N}}

\newcommand{\PP}{\mathbb{P}}

\newcommand{\RR}{\mathbb{R}}

\newcommand{\bone}{\mathbf{1}}

\theoremstyle{plain}
\newtheorem{theorem}{Theorem}[section]

\newtheorem{lemma}[theorem]{Lemma}
\newtheorem{proposition}[theorem]{Proposition}

\newtheorem{ass}[theorem]{Assumption}

\theoremstyle{definition}
\newtheorem{remark}[theorem]{Remark}
\newtheorem{definition}[theorem]{Definition}

\numberwithin{equation}{section}

\begin{document}
\title[Stefan problem with surface tension]{Stefan problem with surface tension: global existence of physical solutions under radial symmetry} 
\author{Sergey Nadtochiy and Mykhaylo Shkolnikov}
\address{Department of Applied Mathematics, Illinois Institute of Technology, Chicago, IL 60616.}
\email{snadtochiy@iit.edu}
\address{ORFE Department, Bendheim Center for Finance, and Program in Applied \& Computational Mathematics, Princeton University, Princeton, NJ 08544.}
\email{mshkolni@gmail.com}
\footnotetext[1]{S.~Nadtochiy is partially supported by the NSF CAREER grant DMS-1651294.}
\footnotetext[2]{M.~Shkolnikov is partially supported by the NSF grant DMS-2108680.}

\begin{abstract}
	We consider the Stefan problem with surface tension, also known as the Stefan-Gibbs-Thomson problem, in an ambient space of arbitrary dimension. Assuming the radial symmetry of the initial data we introduce a novel ``probabilistic'' notion of solution, which can accommodate the discontinuities in time (of the radius) of the evolving aggregate. Our main result establishes the global existence of a probabilistic solution satisfying the natural upper bound on the sizes of the discontinuities. Moreover, we prove that the upper bound is sharp in dimensions $d\ge3$, in the sense that none of the discontinuities in the solution can be decreased in magnitude. The detailed analysis of the discontinuities, via appropriate stochastic representations, differentiates this work from the previous literature on weak solutions to the Stefan problem with surface tension. 
\end{abstract}

\maketitle


\section{Introduction}

Free boundary problems for the heat equation have been introduced independently by \textsc{Lam\'{e}} \& \textsc{Clapeyron} in \cite{LC} and by \textsc{Stefan} in \cite{Stefan1,Stefan2,Stefan3,Stefan4}, and are now commonly known as Stefan problems.~Following a lecture by \textsc{Brillouin} at the Institute Henri Poincar\'{e} in 1929 and its publication \cite{Bri} Stefan problems have taken a central place in the theory of partial differential equations (see, e.g., the classical reference \cite{Rub}) and in mathematical physics (see, e.g., the classical reference \cite{Vis}).~They are now considered canonical as models of melting and freezing (and also of evaporation and condensation).~The original Stefan problem is endowed with a Dirichlet boundary condition that corresponds to melting and freezing at a constant equilibrium temperature.~This is, however, in violation of the Gibbs-Thomson principle which asserts that smaller convex solid crystals find themselves in equilibrium with the associated liquid at a \textit{lower} temperature, due to a surface tension effect (see, e.g., \cite[Subsection 8.3.1]{Gli} for more details).~Stefan problems taking the Gibbs-Thomson principle into account are referred to as \textit{Stefan problems with surface tension} or \textit{Stefan-Gibbs-Thomson problems} (see, e.g., \cite[Problem 1.1]{Vis}). 

\medskip

The classical formulation of the Stefan problem with surface tension can be stated as follows. Given some $d\geq1$, $\Gamma_0\subset\RR^d$, $v(0,\cdot)\!:\RR^d\to\RR$ and $T\in(0,\infty)$, one needs to find $(\Gamma_t\subset\RR^d)_{t\in(0,T)}$ and $(v(t,\cdot)\!:\RR^d\to\RR)_{t\in(0,T)}$ satisfying
\begin{eqnarray}
&& \partial_t v = \frac{1}{2}\Delta v,\quad y\notin\partial\Gamma_t,\;\; t\in(0,T), \label{heat_eq} \\
&& v=H,\quad y\in\partial\Gamma_t,\;\; t\in(0,T), \label{Gibbs-Thomson_eq} \\
&& V=\frac{1}{2}(\nabla v)\cdot\overrightarrow{n_+}+\frac{1}{2}(\nabla v)\cdot\overrightarrow{n_-}
,\quad y\in\partial\Gamma_t,\;\; t\in[0,T), \label{Stefan_eq}
\end{eqnarray}
where $H$ is proportional to the mean curvature of $\partial\Gamma_t$ (chosen so that $H\ge0$ if $\Gamma_t$ is convex), $V$ is the normal growth speed of $\Gamma_t$, and $\overrightarrow{n_+}$, $\overrightarrow{n_-}$ are the normal vector fields on $\partial\Gamma_t$ directed into $\RR^d\backslash\Gamma_t$, $\Gamma_t$, respectively.~In physical terms, for each $t\in[0,T)$, the set $\Gamma_t$ describes the region occupied by a solid, whereas $\RR^d\backslash\Gamma_t$ is the region occupied by the associated liquid, and $v(t,\cdot)$ captures the temperature distribution \textit{below} the equilibrium freezing point for the material in question.~Equation \eqref{heat_eq} then postulates that the temperature distribution evolves according to the (standard) heat equation; equation \eqref{Gibbs-Thomson_eq} quantifies the Gibbs-Thomson principle by enforcing a freezing temperature $H$ degrees \textit{below} the equilibrium freezing point; finally, equation \eqref{Stefan_eq}, known as the Stefan growth condition, asserts that the solid grows (or shrinks) according to the sum of the (negative) temperature slopes in the liquid and in the solid. 

\medskip

A global classical solution to \eqref{heat_eq}--\eqref{Stefan_eq} fails to exist in general, even when the initial data $(\Gamma_0,v(0,\cdot))$ is radially symmetric (see \cite{meir}).~This is due to temperatures far \textit{below} the equilibrium freezing point in regions of high curvature along the solid-liquid interface where melting at an infinite rate may occur.~While the existence of such blow-ups is established in \cite{meir}, the author points out that ``it is difficult to see what happens to the solution'' after a blow-up occurs (see also \cite{GP} for a related discussion in the setting of $d=3$).~This feature has led to work on \eqref{heat_eq}--\eqref{Stefan_eq} in specific situations where well-behaved solutions do exist, e.g., for small surface tension parameters $\gamma>0$ (see \cite{FrRe}), under a smallness assumption on the initial data (see \cite{eps}), near a flat initial interface (see \cite{HaGu}), or when the initial interface is close to a steady sphere (see \cite{ha}).~In the full generality of~\eqref{heat_eq}--\eqref{Stefan_eq} weak solutions, possibly with a phase function taking values strictly between $0$ and $1$ (i.e., a non-sharp interface), have been shown to exist globally in \cite{vis1989}. In the seminal paper \cite{luc}, the global existence of weak solutions with a sharp interface is proved (see also \cite[Example 5]{RoSa} and the references therein for an alternative proof via an abstract gradient flow approach).~As noted in \cite[Section 5]{luc}, the solution concept in \cite{luc} is too weak to yield a unique solution.~We therefore argue the global existence for a new \textit{stronger} notion of a weak solution to \eqref{heat_eq}--\eqref{Stefan_eq}, focusing on initial data that is radially symmetric.

\medskip

The presence of the Gibbs-Thomson condition \eqref{Gibbs-Thomson_eq} makes it impossible to use the global comparison principle that is available for the (radially symmetric) Stefan problem with the Dirichlet boundary condition. This complicates the mathematical analysis of \eqref{heat_eq}--\eqref{Stefan_eq} and, in some sense, is responsible for the formation of singularities in time. On the other hand, the Gibbs-Thomson condition ensures the absence of singularities in space, i.e., the boundary $\partial\Gamma_t$ remains sufficiently regular. The latter is important for establishing the well-posedness of Stefan-type problems in multiple space dimensions: e.g., the results of \cite{nsz} illustrate how the lack of regularity in space may cause an approximating scheme to converge to a wrong limit.~Although the regularity in space is irrelevant in the radially symmetric case (as the boundary of $\Gamma_t$ is always a sphere), it is important to stress that our main motivation for the present work is to develop methods that can ultimately be used to prove the well-posedness of \eqref{heat_eq}--\eqref{Stefan_eq} \textit{without} any spatial symmetry assumptions. In fact, Section~\ref{se:2} provides a ``forward representation" (akin to a particle system) of the proposed Euler scheme for \eqref{heat_eq}--\eqref{Stefan_eq} (see Lemma \ref{le:Sec2.forward}) which can be implemented without the assumption of radial symmetry.~Of course, the proof of convergence for this approximation in the absence of radial symmetry requires additional tools and does not follow immediately from the present work.


\medskip

From here on we assume that $\Gamma_0$ is a ball of some radius $\Lambda_0>0$ and that $v(0,\cdot)$ is radially symmetric.~It is then natural to look for solutions to \eqref{heat_eq}--\eqref{Stefan_eq} for which every $\Gamma_t$ is a ball and all $v(t,\cdot)$ are radially symmetric.~Letting $\Lambda_t$ be the radius of $\Gamma_t$ and taking $u(t,|y|)=v(t,y)$ we recast the problem \eqref{heat_eq}--\eqref{Stefan_eq}, with a minor abuse of notation, as
\begin{eqnarray}
&& \partial_t u = \frac{1}{2}\partial_{xx}u+\frac{d-1}{2x}\,\partial_x u,\quad x\neq\Lambda_t,\;\;t\in(0,T), 
\label{heat_eq_rad} \\
&& u(t,\Lambda_t)=\frac{\gamma}{\Lambda_t}=:H(\Lambda_t),\quad t\in(0,T), \label{Gibbs-Thomson_rad} \\
&& \Lambda'_t=\frac{1}{2}u_x(t,\Lambda_t+)-\frac{1}{2}u_x(t,\Lambda_t-),\quad t\in[0,T). \label{Stefan_rad}
\end{eqnarray}   
Writing $\RR_+$ for $[0,\infty)$ we call a function $\Lambda$ in the Skorokhod space $D([0,T),\RR_+)$ a \textit{probabilistic solution} to \eqref{heat_eq_rad}--\eqref{Stefan_rad} on $[0,T)$ with initial data $(\Lambda_{0-},u(0-,\cdot))$ if 
\begin{equation}\label{what_is_u_intro}
u(t,x):=\E^x\big[\mathbf{1}_{\{\tau_{\Lambda_{t-\cdot}}\leq t\}}\,H(R_{\tau_{\Lambda_{t-\cdot}}})\big]
+\E^x\big[\mathbf{1}_{\{\tau_{\Lambda_{t-\cdot}}>t\}}\,u(0-,R_t)\big],\quad (t,x)\in[0,T)\times\RR_+
\end{equation}
satisfies 
\begin{equation}\label{Stefan_growth_weak}
	\frac{1}{d}\big((\Lambda_t)^d-(\Lambda_{0-})^d\big)
	=\int_{\RR_+} u(0-,x)\,\nu(\mathrm{d}x)-\int_{\RR_+} u(t,x)\,\nu(\mathrm{d}x),\quad t\in[0,T\wedge\zeta), 
\end{equation}
with $R$ being a $d$-dimensional Bessel process started from $x$ under $\PP^x$,
\begin{align}
&\tau_{\Lambda_{t-\cdot}}:=\inf\{s\in[0,t+1]:\,(R_s-\Lambda_{t-s})(x-\Lambda_t)\le 0\},\label{eq.Sec1.tau.def}\\
& \nu(\mathrm{d}x):=x^{d-1}\,\mathrm{d}x,\quad \zeta:= \inf\{t\in[0,T]:\,\Lambda_t = 0\}, \nonumber
\end{align}
and where we set $\Lambda_t:=\Lambda_{0-}$, $t\in[-1,0)$ for convenience.~Hereby, the choice of $R$ mirrors \eqref{heat_eq_rad} (note that, for $d=1$, we implicitly have $u_x(t,0)=0$, since then $R$ is a reflected standard Brownian motion); the definition \eqref{what_is_u_intro} encodes the boundary condition \eqref{Gibbs-Thomson_rad} together with the initial condition $u(0-,\cdot)$; and the growth condition \eqref{Stefan_growth_weak} is a weak formulation of the growth condition \eqref{Stefan_rad}. 

\medskip

To see why \eqref{what_is_u_intro}--\eqref{Stefan_growth_weak} is a natural weak version of \eqref{heat_eq_rad}--\eqref{Stefan_rad}, consider any classical solution $(\Lambda,u)$ to the latter system, such that, e.g., $|u(t,x)|$ decays exponentially in $x$ locally uniformly in $t\in[0,T)$, $x^{d-1}\partial_x u(t,x)$ converges to $0$ as $x\rightarrow\infty$ and as $x\downarrow0$, and $|\partial_t u|$ is $\nu(\mathrm{d}x)\times\mathrm{d}t$ integrable.~Then, \eqref{what_is_u_intro} results from \eqref{heat_eq_rad}, \eqref{Gibbs-Thomson_rad} via the standard Feynman-Kac formula.~To obtain \eqref{Stefan_growth_weak}, we proceed as follows:
\begin{align*}
&\;\int_{\RR_+} u(t,x)\,\nu(\mathrm{d}x)-\int_{\RR_+} u(0,x)\,\nu(\mathrm{d}x)
=\int_0^t\frac{\mathrm{d}}{\mathrm{d}s} \int_0^\infty u(s,x)\,\nu(\mathrm{d}x)\,\mathrm{d}s \\ 
&= \frac{1}{2} \int_0^t \int_0^{\Lambda_s} \left[x^{d-1} \partial_{xx}u+(d\!-\!1)x^{d-2}\,\partial_x u\right]\, \mathrm{d}x\,\mathrm{d}s
+ \frac{1}{2} \int_0^t \int_{\Lambda_s}^{\infty} \left[x^{d-1} \partial_{xx}u+(d\!-\!1)x^{d-2}\,\partial_x u\right]\,\mathrm{d}x\,\mathrm{d}s \\
&= \frac{1}{2} \int_0^t \int_0^{\Lambda_s} \partial_x\left[x^{d-1} \partial_{x}u\right]\,\mathrm{d}x\,\mathrm{d}s
+ \frac{1}{2} \int_0^t \int_{\Lambda_s}^\infty \partial_x\left[x^{d-1} \partial_{x}u\right]\,\mathrm{d}x\,\mathrm{d}s \\ 
&= \frac{1}{2} \int_0^t \Lambda_s^{d-1} \big[\partial_{x}u(s,\Lambda_s-) - \partial_{x}u(s,\Lambda_s+)\big]\,\mathrm{d}s
= - \int_0^t \Lambda_s^{d-1}\,\mathrm{d}\Lambda_s = \frac{1}{d}\big((\Lambda_0)^d - (\Lambda_t)^d\big).
\end{align*}

\medskip

\begin{figure}[h]
\includegraphics[width=8.8cm,height=6.6cm]{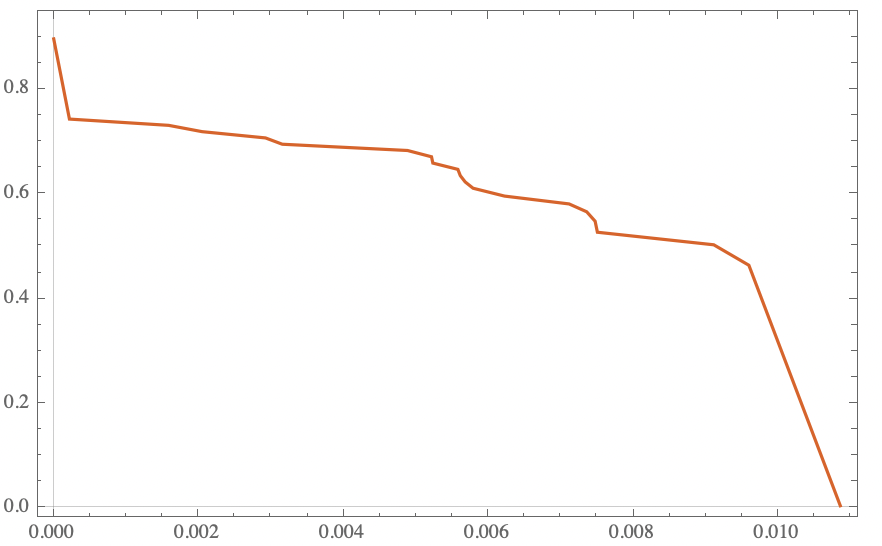}
\caption{The plot depicts a version of the Euler scheme in Definition~\ref{def:Euler}, with Exponential($1/\Delta$) time steps, for $d\!=\!3$, initial data $(\Lambda_{0-},u(0-,\cdot))\!=\!(0.9,\mathbf{1}_{[0,0.81]})$ and parameter $\gamma\!=\!1$.~Hereby, we have chosen $\Delta\!=\!0.5\cdot10^{-3}$ and have discretized space according to a mesh size of $3\cdot10^{-3}$.~The~exponential time steps enable us to use \cite[Section II.5, displays (1.0.4), (1.1.2), (1.1.6), (1.2.2), (1.2.6)]{BorSal} in the implementation.}
\label{fig:two jumps}
\end{figure}

In general, probabilistic solutions $\Lambda$ to \eqref{heat_eq_rad}--\eqref{Stefan_rad} exhibit jumps.~(It is worth mentioning that, unlike the case without surface tension, one cannot exclude the jumps of $\Lambda$ globally even by imposing smallness assumptions on the initial data.)~Figure~\ref{fig:two jumps} shows a numerical simulation of a probabilistic solution $\Lambda$ (via a version of the Euler scheme in Definition~\ref{def:Euler}), which, in particular, suggests multiple downward jumps and illustrates the complexity of the dynamics of $\Lambda$, even for very simple initial data.~Note also that a jump may occur at time zero, which explains the notation $(\Lambda_{0-},u(0-,\cdot))$ for the input data.~(The energy of the system is preserved at the initial jump, via \eqref{Stefan_growth_weak}.) At any given time, many jump sizes may be consistent with the growth condition~\eqref{Stefan_growth_weak}.~The following additional condition rules out the ambiguity about the jump sizes in $\Lambda$: 
\begin{align}
&\Lambda_{t-}-\Lambda_t =  \inf\bigg\{y\in(0,\Lambda_{t-}]:\, \int_{\Lambda_{t-}-y}^{\Lambda_{t-}} u(t-,x)\,\nu(\mathrm{d}x) >\int_{\Lambda_{t-}-y}^{\Lambda_{t-}} \big(H(x)-1\big)\,\nu(\mathrm{d}x) \bigg\},\label{phys_jump}\\
& \Lambda_{t}-\Lambda_{t-} = \inf\bigg\{y>0:\, \int_{\Lambda_{t-}}^{\Lambda_{t-}+y} u(t-,x)\,\nu(\mathrm{d}x) <\int_{\Lambda_{t-}}^{\Lambda_{t-}+y} \big(H(x)+1\big)\,\nu(\mathrm{d}x)\bigg\}.\label{phys_jump.n}
\end{align}
We call a probabilistic solution $\Lambda$ \textit{physical} if it satisfies \eqref{phys_jump} at any downward jump time~$t$ and \eqref{phys_jump.n} at any upward jump time $t$.~It is shown in Subsections \ref{subse:Sec4.prelim} and \ref{subse:lowerBound} that, under the Assumption \ref{main_ass} below and under a natural upper bound on the absolute jump sizes, any probabilistic solution $\Lambda$ has no upward jumps and that the absolute sizes of its downward jumps are bounded from below by the right-hand side of \eqref{phys_jump}, provided $d\geq3$.~Thus, a physical solution has the smallest possible jumps (which correspond to the smallest amounts of energy exchange between the phase configuration and the temperature).
~The physicality condition \eqref{phys_jump}--\eqref{phys_jump.n} is analogous with the one recently leveraged to uniquely determine the solution of the one-dimensional one-phase supercooled Stefan problem without surface tension (see \cite[display (1.5)]{dns}, as well as \cite{DIRT1}, \cite{NaSh1}, \cite{HLS}, \cite{LS} that led up to \cite{dns}). 

\medskip

We work throughout under the following assumption on the input data $(\Lambda_{0-},u(0-,\cdot))$. 

\begin{ass}\label{main_ass}
The function $u(0-,\cdot)$ is non-negative, continuous, $\int_{\RR_+} u(0-,x)\,\nu(\mathrm{d}x)\!=\!1$, and there exist constants $C_1\in(0,1)$, $C_2>0$, $\alpha>0$ such that $u(0-,x)\le C_1\wedge \left(C_2e^{-\alpha x}\right)$, $x\ge0$. In addition, $\Lambda_{0-}>0$.
\end{ass}

The above assumption holds throughout the paper even if not cited explicitly.
We are now ready to state our main result. 

\begin{theorem}\label{thm_main}
For any $\gamma>0$, any $(\Lambda_{0-},u(0-,\cdot))$ satisfying Assumption \ref{main_ass}, and any $T\in(0,\infty)$, there exists a probabilistic solution $\Lambda\in D([0,T),\RR_+)$, in the sense that $\Lambda$ satisfies \eqref{Stefan_growth_weak} for all $t\in[0,T\wedge\zeta)$, with $u$ given by \eqref{what_is_u_intro}, and it is such that $\Lambda$ has no upward jumps and \eqref{phys_jump} holds with ``$=$'' replaced by ``$\leq$'' for all $t\in[0,T\wedge\zeta]$.
If, in addition, $d\ge3$, then $\Lambda$ can be chosen so that it is physical, in the sense that \eqref{phys_jump} holds for all $t\in[0,T\wedge\zeta)$.
\end{theorem}

\begin{remark}
A few comments on Assumption \ref{main_ass} and on Theorem \ref{thm_main} are in order:
\begin{enumerate}[(a)]
\item The non-negativity of $u(0-,\cdot)$ means that the liquid is \textit{supercooled} globally, and not only near the solid-liquid interface as required by the Gibbs-Thomson condition \eqref{Gibbs-Thomson_rad}.~This assumption can be easily relaxed to allow for $u(0-,\cdot)$ that also takes negative values but is bounded from below by a constant $-a$. This relaxation is achieved by a simple change of variables $u\mapsto u+a$ and by setting $H(x)=\frac{\gamma}{x} + a$.
\item The assumption $\int_{\RR_+} u(0-,x)\,\nu(\mathrm{d}x)=1$ can be relaxed to $\int_{\RR_+} u(0-,x)\,\nu(\mathrm{d}x)\in(0,\infty)$ by a change of the space scale. 
\item The inequality $u(0-,x)\le C_1\wedge \left(C_2e^{-\alpha x}\right)$ can be broken down into $\|u(0-,\cdot)\|_\infty<1$ and an exponential decay of $u(0-,\cdot)$ at infinity.~The former excludes the so-called \textit{hypercooled} regime, in which distinct phenomena are observed experimentally, such as glass formation (see \cite[Subsection 17.3.2]{Gli} for more details). 
\item The problem can be set up on a bounded domain (as opposed to $\RR_+$) by prescribing a Dirichlet or a Neumann condition at the exterior boundary (which would correspond to absorption or reflection of the associated Bessel process).
\item The proof of Theorem \ref{thm_main} is constructive and uses a family of numerical approximations whose limit points are probabilistic solutions.
\item Note that the size of a jump at the time $\zeta$ of complete melting may not satisfy \eqref{phys_jump}, and the growth (or, energy preservation) condition \eqref{Stefan_growth_weak} may fail at $t=\zeta$.~This is because the radius of the ball at $\zeta-$ may be too small to realize the prescribed jump and to fulfill the growth condition.~(Physically, this means that the temperature along $[0,\Lambda_{\zeta-}]$ will not be equal to $H$.)~However, it is important to notice that the absolute jump size at $\zeta$ is still bounded from above by the right-hand side of \eqref{phys_jump}, which excludes ``unnecessary" jumps and, in turn, the obvious cases of non-uniqueness.
\item We conjecture that the additional physicality condition \eqref{phys_jump} leads to uniqueness of a probabilistic solution $\Lambda$.~This conjecture is a subject of ongoing research.
\item We also believe that the assumption $d\geq3$ can be relaxed, but this relaxation requires the use of somewhat different methods and is also a subject of ongoing investigation.
\end{enumerate}
\end{remark}

The rest of the paper is structured as follows.~In Section \ref{se:2}, we introduce (implicit) Euler scheme approximations associated with \eqref{what_is_u_intro}--\eqref{Stefan_growth_weak}.~We then show the relative compactness of the Euler scheme approximations with respect to the Skorokhod M1 topology, by passing to their equivalent \textit{forward} probabilistic formulation.~A particular challenge is the non-monotonicity of~$\Lambda$, in contrast to the setting of the one-dimensional one-phase supercooled Stefan problem studied in \cite{dns}.~In Section \ref{se:3}, we investigate the limits of the Euler scheme approximations and derive the weak Stefan growth condition \eqref{Stefan_growth_weak}.~Hereby, we use the backward probabilistic formulation together with the remarkable observation of \cite[Proposition 3.4(i)]{ABPR}.~In Section~\ref{se:4}, we verify the physicality condition \eqref{phys_jump}: in Subsections~\ref{subse:Sec4.prelim} and \ref{subse:lowerBound} we prove that every solution to \eqref{what_is_u_intro}--\eqref{Stefan_growth_weak} satisfies the desired lower bound on the absolute sizes of its downward jumps, provided $d\geq3$, and Subsection \ref{subse:upperBound} shows that the jumps of the limit points of the proposed Euler approximations satisfy the same upper bound and that they do not jump upwards.~The former is proved by suitably recasting the right-hand side in \eqref{what_is_u_intro} and by relying on the supermartingale property of $H(R)$ when $d\ge3$.~The latter is obtained by establishing a corresponding upper bound on the time increments of the Euler scheme approximations and that this upper bound is preserved in the limit.

\section{Euler scheme approximation: definition and relative compactness} \label{se:2}

Our starting point is an (implicit) Euler scheme approximation for \eqref{what_is_u_intro}--\eqref{Stefan_growth_weak}.

\begin{definition}\label{def:Euler}
\!\! For $\Delta\!>\!0$, set $\Lambda^\Delta_0\!=\!\Lambda_{0-}$, $u^\Delta(0,x)\!=\!u(0-,x)$ and, for $m\!=\!0,1,\ldots,\lceil T/\Delta\rceil-1$,
\begin{eqnarray*}
&&
u^\Delta(t,x)=\E^x\big[\mathbf{1}_{\{\tau_{\Lambda^\Delta_{m\Delta}}<t-m\Delta\}}
\,H(\Lambda^\Delta_{m\Delta})\big]
+\E^x\big[\mathbf{1}_{\{\tau_{\Lambda^\Delta_{m\Delta}}\ge t-m\Delta\}}\,u^\Delta(m\Delta,R_{t-m\Delta})\big], \\
&& \qquad \qquad \qquad \qquad \qquad \qquad \qquad \qquad \qquad \qquad \qquad \qquad \;\, t\in(m\Delta,(m+1)\Delta),\,\,x>0, \\ 
&&\,\,\\
&& \Lambda^\Delta_{(m+1)\Delta}=\begin{cases}
0\quad\text{if}\;\;\;\Lambda^\Delta_{m\Delta}=0, \\
\,\,\\
0\vee\sup\big\{y\in[0,\Lambda^\Delta_{m\Delta})\!: \\
\;\,\int_{\RR_+}\! u^\Delta((m\!+\!1)\Delta-,x)\,\mathbf{1}_{\RR_+\backslash [y,\Lambda^\Delta_{m\Delta}]}(x)
\!+\! H(x)\,\mathbf{1}_{[y,\Lambda^\Delta_{m\Delta}]}(x)
\!-\! u^\Delta(m\Delta,x)\,\nu(\mathrm{d}x) \\ 
\qquad\qquad\qquad\qquad\qquad\qquad\qquad\qquad\qquad\qquad\qquad\qquad\;\;
<\!\frac{1}{d}\big((\Lambda^\Delta_{m\Delta})^d\!-\!y^d\big)\big\} \\
\;\;\;
\text{if}\;\;\;\Lambda^\Delta_{m\Delta}>0\;\;\;\text{and}\;\;\;\int_{\RR_+} u^\Delta((m+1)\Delta-,x)\,\nu(\mathrm{d}x)\ge
\int_{\RR_+} u^\Delta(m\Delta,x)\,\nu(\mathrm{d}x),  \\
\,\,\\
\inf\big\{y>\Lambda^\Delta_{m\Delta}\!: \\
\;
\int_{\RR_+} u^\Delta((m\!+\!1)\Delta-,x)\,\mathbf{1}_{\RR_+\backslash[\Lambda^\Delta_{m\Delta},y]}(x)
\!+\! H(x)\,\mathbf{1}_{[\Lambda^\Delta_{m\Delta},y]}(x)
\!-\! u^\Delta(m\Delta,x)\,\nu(\mathrm{d}x) \\ 
\qquad\qquad\qquad\qquad\qquad\qquad\qquad\qquad\qquad\qquad\qquad\qquad\;\;
>\!\frac{1}{d}\big((\Lambda^\Delta_{m\Delta})^d\!-\!y^d\big)\big\} \\
\;\;\;
\text{if}\;\;\;\Lambda^\Delta_{m\Delta}>0\;\;\;\text{and}\;\;\;\int_{\RR_+} u^\Delta((m+1)\Delta-,x)\,\nu(\mathrm{d}x)<
\int_{\RR_+} u^\Delta(m\Delta,x)\,\nu(\mathrm{d}x),
\end{cases} \\
&&\,\,\\
&& u^\Delta((m+1)\Delta,x)=\begin{cases}
u^\Delta((m+1)\Delta-,x)\,\mathbf{1}_{\RR_+\backslash[\Lambda^\Delta_{(m+1)\Delta},\Lambda^\Delta_{m\Delta}]}(x)
+ H(x)\,\mathbf{1}_{[\Lambda^\Delta_{(m+1)\Delta},\Lambda^\Delta_{m\Delta}]}(x)\\
\qquad\qquad\qquad\qquad\qquad\qquad\qquad\qquad\qquad\quad
\text{if}\;\;\;\Lambda^\Delta_{(m+1)\Delta}\le \Lambda^\Delta_{m\Delta}, \\
\,\,\\
u^\Delta((m+1)\Delta-,x)\,\mathbf{1}_{\RR_+\backslash[\Lambda^\Delta_{m\Delta},\Lambda^\Delta_{(m+1)\Delta}]}(x)
+ H(x)\,\mathbf{1}_{[\Lambda^\Delta_{m\Delta},\Lambda^\Delta_{(m+1)\Delta}]}(x)\\
\qquad\qquad\qquad\qquad\qquad\qquad\qquad\qquad\qquad\quad
\text{if}\;\;\;\Lambda^\Delta_{(m+1)\Delta}>\Lambda^\Delta_{m\Delta}.
\end{cases}
\end{eqnarray*}
In the above, $R$ is a Bessel process of dimension $d$ started from $x\ge0$ under $\PP^x$, $\tau_\ell$ is its first hitting time of level $\ell>0$, and $\nu(\mathrm{d}x):=x^{d-1}\,\mathrm{d}x$ on $\RR_+$.~We also use the following conventions: $\tau_0:=\infty$, $H(0):=0$ (only needed for $d=1$), $\sup\emptyset=-\infty$, $\inf\emptyset=\infty$.

Finally, we define $\Lambda^\Delta$ as the unique extension of the above values $\Lambda^\Delta_{m\Delta}$ to a right-continuous function on $[-1,T+1]$ that is constant on each interval $(m\Delta,(m+1)\Delta)\cap(0,T)$, and is equal to $\Lambda_{0-}$ in $[-1,0)$ and to $\Lambda^\Delta_T$ in $[T,T+1]$.
\end{definition}

\begin{remark}
In the expressions appearing in Definition \ref{def:Euler} we implicitly rely on the fact that $u^\Delta((m+1)\Delta-,x)=\lim_{s\uparrow (m+1)\Delta}u^{\Delta}(s,x)$ is well-defined. The proof of this fact can be carried out in the same way as the proof of Lemma \ref{le:u.LRlims} below and is omitted here for brevity.
\end{remark}

The next proposition is the main result of this section. 

\begin{proposition}\label{prop:tightness}
The family $\{\Lambda^\Delta\}_{\Delta\downarrow0}$ is relatively compact, in the sense that every sequence $(\Lambda^{\Delta_n})_{n\ge1}$ with $\Delta_n\downarrow0$ has a subsequence converging in the Skorokhod space $D([-1,T+1],\RR_+)$ equipped with the M1 topology.
\end{proposition}

The proof of Proposition \ref{prop:tightness} relies on the following two lemmas, which invoke the next probabilistic objects: 
\begin{itemize}
\item a $d$-dimensional Bessel process $X$ started according to the density $u(0-,\cdot\,)\,x^{d-1}$; 
\item its hitting time $\tau^\Delta:=\inf\{t>0\!:(X_t-\Lambda^\Delta_t)(X_0-\Lambda^\Delta_0)<0\}$;
\item for every $m\ge 1$, the $d$-dimensional Bessel processes $\{X^{m,i,\Delta}\}_{i\ge1}$ started at time $m\Delta$ from the atoms of an independent Poisson random measure of intensity $H$ in the interval between $\Lambda^\Delta_{(m-1)\Delta}$ and $\Lambda^\Delta_{m\Delta}$;   
\item their hitting times $\tau^{m,\Delta}_i:=\inf\{t>m\Delta\!:(X^{m,i,\Delta}_t-\Lambda^\Delta_t)(X^{m,i,\Delta}_{m\Delta}-\Lambda^\Delta_{m\Delta})<0\}$;
\item the jumps times $\{T^{\delta,\Delta}_i\}_{i\ge 1}$ of a Poisson process with rate $2\gamma\delta^{-1}(\Lambda^\Delta)^{d-2}$, for $\delta>0$; 
\item $[-1,1]$-valued independent uniform random variables $\{\gamma_i\}_{i\ge 1}$; 
\item independent $d$-dimensional Bessel processes $\{Y^{\delta,i,\Delta}\}_{i\ge1}$ started at the times $\{T^{\delta,\Delta}_i\}_{i\ge1}$ from $\{(\Lambda^\Delta_{T^{\delta,\Delta}_i}+\delta\gamma_i)\vee0\}_{i\ge 1}$, respectively;
\item their hitting times $\tau^{\delta,\Delta}_i:=\inf\{t>T^{\delta,\Delta}_i\!:(Y^{\delta,i,\Delta}_t-\Lambda^\Delta_t)(Y^{\delta,i,\Delta}_{T^{\delta,\Delta}_i}-\Lambda^\Delta_{T^{\delta,\Delta}_i})<0\}$. 
\end{itemize}

\begin{lemma}\label{le:Sec2.forward}
Fix a $\Delta>0$. Then, for all $m=0,\,1,\,\ldots,\,\inf\{m_0\ge1\!:\Lambda^\Delta_{m_0\Delta}=0\}-1$, 
\begin{equation}
\begin{split}
& \frac{1}{d}\big((\Lambda^\Delta_{m\Delta})^d\!-\!(\Lambda^\Delta_{0-})^d\big)
=\PP(\tau^\Delta\!\le\! m\Delta)
-\sum_{n=1}^m\,\sum_{i\ge 1}\, \PP(\tau^{n,\Delta}_i\!>\! m\Delta)
-\lim_{\delta\downarrow0}\,\sum_{i\ge1}\,\PP(T^{\delta,\Delta}_i\!\le\! m\Delta\!<\!\tau^{\delta,\Delta}_i),
\label{Lambda_forward} \\
& u^\Delta(t,\cdot)\,\nu\!=\!\PP(X_t\!\in\!\cdot\,,\tau^\Delta\!>\!t)
\!+\!\sum_{n=1}^{\lfloor t/\Delta\rfloor}\!
\sum_{i\ge1} \PP(X^{n,i,\Delta}_t\!\in\!\cdot\,,\tau^{n,\Delta}_i\!>\!t)
\!+\!\lim_{\delta\downarrow0}\sum_{i\ge1}\PP(Y^{\delta,i,\Delta}_t\!\in\!\cdot\,,T_i^{\delta,\Delta}\!\le\! t\!<\!\tau^{\delta,\Delta}_i), \\
&\qquad\qquad\qquad\qquad\qquad\qquad\qquad\qquad\qquad\qquad\qquad\qquad\qquad\qquad\qquad\,
t\in[(m\!-\!1)^+\Delta,m\Delta],
\end{split}
\end{equation}
where the second $\delta\downarrow0$ limit is in the sense of the total variation convergence of measures.

\end{lemma}

\noindent\textbf{Proof.} We argue \eqref{Lambda_forward} by induction over $m$. The two assertions are clear for $m\!=\!0$. Suppose now that we have \eqref{Lambda_forward} with $(m-1)$ in place of $m$, for an $m\ge1$. Let us show that the first line of \eqref{Lambda_forward} holds and that
\begin{equation}\label{forward_to_backward}
\begin{split}
u^\Delta(t,\cdot)\,\nu\!=\!\PP(X_t\!\in\!\cdot\,,\tau^\Delta\!>\!t)
\!+\!\sum_{n=1}^{\lfloor t/\Delta\rfloor}\!
\sum_{i\ge1} \PP(X^{n,i,\Delta}_t\!\in\!\cdot\,,\tau^{n,\Delta}_i\!>\!t)
\!+\!\lim_{\delta\downarrow0}\sum_{i\ge1}\PP(Y^{\delta,i,\Delta}_t\!\in\!\cdot\,,T_i^{\delta,\Delta}\!\le\! t\!<\!\tau^{\delta,\Delta}_i), \\ 
t\in((m-1)\Delta,m\Delta],
\end{split}
\end{equation}
where the $\delta\downarrow0$ limit is in the sense of the total variation convergence of measures.

\smallskip

\noindent\textbf{Step 1.}
We set $\lambda=\Lambda^\Delta_{(m-1)\Delta}$ and use the induction hypothesis along with the Markov property of the Bessel process to compute, for any $s\in(0,\Delta)$, any Borel $A\subset\RR_+$, and $t:=(m-1)\Delta+s$,
\begin{equation*}
\begin{split}
& \;\PP(X_t\!\in\! A,\tau^\Delta\!>\!t)
\!+\!\sum_{n=1}^{m-1} \! \sum_{i\ge1} \PP(X^{n,i,\Delta}_t\!\in\!A,\tau^{n,\Delta}_i\!>\!t)
\!+\!\lim_{\delta\downarrow0}\sum_{i\ge1} \PP(Y^{\delta,i,\Delta}_t\!\in\! A,
T_i^{\delta,\Delta}\!\le\!(m\!-\!1)\Delta,t\!<\!\tau^{\delta,\Delta}_i) \\
&=\int_{\RR_+} \PP^x(R_s\in A,\,\tau_\lambda> s) \Big[\PP\big(X_{(m-1)\Delta}\!\in\! \mathrm{d}x,\,\tau^\Delta\!>\!(m-1)\Delta\big) \\
&\quad\!+\!\sum_{n=1}^{m-1} \! \sum_{i\ge1} \PP\big(X^{n,i,\Delta}_{(m-1)\Delta}\!\in\!\mathrm{d}x,\tau^{n,\Delta}_i\!>\!(m\!-\!1)\Delta\big) 
\!+\!\lim_{\delta\downarrow0}\sum_{i\ge1} \PP\big(Y^{\delta,i,\Delta}_{(m-1)\Delta}\!\in\!\mathrm{d}x,
T_i^{\delta,\Delta}\!\le\!(m\!-\!1)\Delta<\!\tau^{\delta,\Delta}_i\big)\Big]\\
& = \int_{\RR_+} \PP^x(R_s\in A,\,\tau_\lambda> s)\,u^\Delta((m-1)\Delta,x)\,\nu(\mathrm{d}x) \\
& = \int_{\RR_+} \int_A \psi(s;x,y)\,
\PP^{x\to y}(\tau_\lambda>s)\,\mathrm{d}y\,u^\Delta((m-1)\Delta,x)\,\nu(\mathrm{d}x),
\end{split}
\end{equation*}
where $\psi$ is the transition kernel of the $d$-dimensional Bessel process and $\PP^{x\to y}$ is the law of the $d$-dimensional Bessel bridge from $x$ to $y$ on $[0,s]$. The explicit formula for $\psi$ (see, e.g., \cite[display (17)]{Law}) shows that $x^{d-1}\,\psi(s;x,y)=y^{d-1}\,\psi(s;y,x)$ and that under $\PP^{x\to y}$ the time-reversal of the canonical process has the law $\PP^{y\to x}$. These and Fubini's Theorem yield
\begin{equation*}
\begin{split}
& \;\PP(X_t\!\in\! A,\tau^\Delta\!>\!t)
\!+\!\sum_{n=1}^{m-1} \! \sum_{i\ge1} \PP(X^{n,i,\Delta}_t\!\in\!A,\tau^{n,\Delta}_i\!>\!t)
\!+\!\lim_{\delta\downarrow0}\sum_{i\ge1} \PP(Y^{\delta,i,\Delta}_t\!\in\! A,
T_i^{\delta,\Delta}\!\le\!(m\!-\!1)\Delta,t\!<\!\tau^{\delta,\Delta}_i) \\
&=\, \int_A \int_{\RR_+} \psi(s;y,x)\,\PP^{y\to x}(\tau_\lambda>s)\,u^\Delta((m-1)\Delta,x)\,\mathrm{d}x\,\nu(\mathrm{d}y) \\
&=\, \int_A \E^y\big[\mathbf{1}_{\{\tau_\lambda>s\}}\,u^\Delta((m-1)\Delta,R_s)\big]\,\nu(\mathrm{d}y).
\end{split}
\end{equation*}

\smallskip

\noindent\textbf{Step 2.} Next, we notice that
\begin{equation}\label{eq.Sec2.mainLemma.Step2.eq1}
\begin{split}
&\;\sum_{i\ge1}\,\PP(Y^{\delta,i,\Delta}_t\in\cdot\,,\,(m-1)\Delta<T_i^{\delta,\Delta}\le t<\tau^{\delta,\Delta}_i) \\
&=\,\sum_{i\ge1}\,\int_0^s \frac{r^{i-1}\,2^i\,\gamma^i\,\lambda^{(d-2)i}}{(i-1)!\,\delta^i}\,e^{-2\gamma\lambda^{d-2}r/\delta}\,\frac{1}{2\delta}\,\int_{-\delta}^\delta \PP^{\lambda+a}(R_{s-r}\in\cdot\,,\,\tau_\lambda>s-r)\,\mathrm{d}a\,\mathrm{d}r \\
&=\,\frac{\gamma\lambda^{d-2}}{\delta^2}\,\int_0^s \int_{-\delta}^\delta \PP^{\lambda+a}(R_{s-r}\in\cdot\,,\,\tau_\lambda>s-r)\,\mathrm{d}a\,\mathrm{d}r.
\end{split}
\end{equation}
Then, using the notation of Step 1, we conclude that the density of the left-hand side in the above display converges, as $\delta\downarrow0$, to
\begin{equation*}
\begin{split}
&\;\lim_{\delta\downarrow0}\,\frac{\gamma\lambda^{d-2}}{\delta^2}\,\int_0^s \int_{-\delta}^\delta \PP^{\lambda+a\to x}(\tau_\lambda>s-r)\,\psi(s-r;\lambda+a,x)\,\mathrm{d}a\,\mathrm{d}r \\
&=\frac{\gamma x^{d-1}}{\lambda}\,\lim_{\delta\downarrow0} \,\frac{1}{\delta^2}\,\int_0^s \int_{-\delta}^\delta \PP^{x\to\lambda+a}(\tau_\lambda>s-r)\,
\psi(s-r;x,\lambda+a)\,\mathrm{d}a\,\mathrm{d}r \\
&=\frac{\gamma x^{d-1}}{\lambda}\,\lim_{\delta\downarrow0} \,\frac{1}{\delta^2}\,\int_0^s \PP^x(R_{s-r}\in[\lambda-\delta,\lambda+\delta],\,\tau_\lambda>s-r)\,\mathrm{d}r  \\
&=\frac{\gamma x^{d-1}}{\lambda}\,\lim_{\delta\downarrow0} \,\frac{1}{\delta^2}\;
\E^x\bigg[\int_0^{\tau_\lambda\wedge s} \mathbf{1}_{\{R_r\in [\lambda-\delta,\lambda+\delta]\}}\,\mathrm{d}r\bigg],
\end{split}
\end{equation*}
provided the latter limit is well-defined.

\medskip

In view of the occupation time formula (see, e.g., \cite[Chapter 3, Theorem 7.1]{KaSh}) and Fubini's Theorem, this equals to
\begin{equation}\label{loc_time_rep}
\begin{split}
\frac{\gamma x^{d-1}}{\lambda}\,\lim_{\delta\downarrow0} \,\frac{1}{\delta^2}\, \int_{\lambda-\delta}^{\lambda+\delta} 
\E^x\big[L^a_{\tau_\lambda}\,\mathbf{1}_{\{\tau_\lambda\le s\}}\big]\,\mathrm{d}a 
 +  \frac{\gamma x^{d-1}}{\lambda}\,\lim_{\delta\downarrow0} \,\E^x\bigg[\frac{1}{\delta^2}\, \int_{\lambda-\delta}^{\lambda+\delta} 
L^a_s\,\mathbf{1}_{\{\tau_\lambda>s\}}\,\mathrm{d}a\bigg],
\end{split}
\end{equation}
where $L$ is the semimartingale local time of $R$. 

\medskip

For $x>\lambda$, to evaluate the first summand in \eqref{loc_time_rep} we compute, for $d\geq2$ and $\theta>0$:
\begin{equation}\label{lapl_trans_iden}
\begin{split}
&\lim_{\delta\downarrow0} \frac{1}{\delta^2} \int_\lambda^{\lambda+\delta} \E^x\big[L^a_{\tau_\lambda}\,e^{-\theta\tau_\lambda} \big]\,\mathrm{d}a
=\frac{1}{2}\lim_{a\downarrow\lambda}\,\frac{\mathrm{d}}{\mathrm{d}a}\, \E^x\big[L^a_{\tau_\lambda}\,e^{-\theta\tau_\lambda} \big]
= -\frac{1}{2}\lim_{a\downarrow\lambda}\,\lim_{b\downarrow0}
\, \frac{1}{b}\,\frac{\mathrm{d}}{\mathrm{d}a}\, \E^x\big[e^{-\theta\tau_\lambda - b L^a_{\tau_\lambda}} \big]\\
&\!= -\frac{1}{2}\lim_{a\downarrow\lambda}\,\lim_{b\downarrow0}
\, \frac{1}{b}\,\frac{\mathrm{d}}{\mathrm{d}a}\,
\frac{x^{-d/2+1}\,K_{d/2-1}(x\sqrt{2\theta})}{\lambda^{-d/2+1}\,K_{d/2-1}(\lambda\sqrt{2\theta}) 
\!+\! 2b a^{d/2} (2\theta)^{d/2-1} S_{d/2-1}(a\sqrt{2\theta},\lambda\sqrt{2\theta})\, K_{d/2-1}(a\sqrt{2\theta})}\\
&\!= \frac{x^{-d/2+1}\,K_{d/2-1}(x\sqrt{2\theta})}{\lambda^{-d+2}\,K_{d/2-1}(\lambda\sqrt{2\theta})^2}\,(2\theta)^{d/2-1} \,\lim_{a\downarrow\lambda}\,
\frac{\mathrm{d}}{\mathrm{d}a}\,\big[ a^{d/2} S_{d/2-1}(a\sqrt{2\theta},\lambda\sqrt{2\theta})\, K_{d/2-1}(a\sqrt{2\theta}) \big] \\
&\!= \frac{x^{-d/2+1}\,K_{d/2-1}(x\sqrt{2\theta})}{\lambda^{-d+2}\,K_{d/2-1}(\lambda\sqrt{2\theta})^2}\,(2\theta)^{d/2-1} 
\,\lambda^{d/2}\,K_{d/2-1}(\lambda\sqrt{2\theta})\,\lim_{a\downarrow\lambda}\,
\frac{\mathrm{d}}{\mathrm{d}a} S_{d/2-1}(a\sqrt{2\theta},\lambda\sqrt{2\theta})\\
&\!= \frac{x^{-d/2+1}\,K_{d/2-1}(x\sqrt{2\theta})}{\lambda^{-3d/2+2}\,K_{d/2-1}(\lambda\sqrt{2\theta})}\,(2\theta)^{d/2-1} 
\,\sqrt{2\theta}\,\partial_x S_{d/2-1}(\lambda\sqrt{2\theta},\lambda\sqrt{2\theta})\\
&\!= \frac{x^{-d/2+1}\,K_{d/2-1}(x\sqrt{2\theta})}{\lambda^{-3d/2+2}\,K_{d/2-1}(\lambda\sqrt{2\theta})}\,(2\theta)^{d/2-1/2} 
\,(\lambda\sqrt{2\theta})^{1-d}
 = \frac{x^{-d/2+1}\,K_{d/2-1}(x\sqrt{2\theta})}{\lambda^{-d/2+1}\,K_{d/2-1}(\lambda\sqrt{2\theta})}
= \E^x[e^{-\theta\tau_\lambda}],
\end{split}
\end{equation}
where $K_{d/2-1}$ is the modified Bessel function of the second kind, the function $S_{d/2-1}$ is defined in \cite[Appendix 2.10]{BorSal}, the third equality is from \cite[Sections II.4, II.6:~displays~2.3.3]{BorSal}, the last equality is a consequence of \cite[Sections II.4, II.6: displays 2.0.1]{BorSal}, and the intermediate equalities follow from the properties of $S_{d/2-1}$ stated in \cite[Appendix 2.10]{BorSal}.

Similarly, for $d=1$ and $\theta>0$, we use \cite[Section II.3: displays 2.3.3, 2.0.1]{BorSal} to obtain
\begin{equation*}
\begin{split}
&\,\lim_{\delta\downarrow0} \frac{1}{\delta^2} \int_\lambda^{\lambda+\delta} \E^x\big[L^a_{\tau_\lambda}\,e^{-\theta\tau_\lambda} \big]\,\mathrm{d}a
=\frac{1}{2}\lim_{a\downarrow\lambda}\,\frac{\mathrm{d}}{\mathrm{d}a}\, \E^x\big[L^a_{\tau_\lambda}\,e^{-\theta\tau_\lambda} \big]
= -\frac{1}{2}\lim_{a\downarrow\lambda}\,\lim_{b\downarrow0}
\, \frac{1}{b}\,\frac{\mathrm{d}}{\mathrm{d}a}\, \E^x\big[e^{-\theta\tau_\lambda - b L^a_{\tau_\lambda}} \big]\\
&= -\frac{1}{2}\lim_{a\downarrow\lambda}\,\lim_{b\downarrow0}
\, \frac{1}{b}\,\frac{\mathrm{d}}{\mathrm{d}a}\,
\frac{\sqrt{2\theta} \,e^{-(x-a)\sqrt{2\theta}}}{\sqrt{2\theta}\,e^{(a-\lambda)\sqrt{2\theta}} + 2b\,\sinh((a-\lambda)\sqrt{2\theta})}\\
&= -\,\lim_{a\downarrow\lambda}\,\lim_{b\downarrow0}
\, \frac{\sqrt{2\theta}\,e^{-(x-a)\sqrt{2\theta}}}{\sqrt{2\theta}\,e^{(a-\lambda)\sqrt{2\theta}} \!+\! 2b\,\sinh((a\!-\!\lambda)\sqrt{2\theta})}
\frac{\sqrt{2\theta}\,\left[\sinh((a\!-\!\lambda)\sqrt{2\theta})\!-\! \cosh((a\!-\!\lambda)\sqrt{2\theta})\right]}{\sqrt{2\theta}\,e^{(a-\lambda)\sqrt{2\theta}} \!+\! 2b\,\sinh((a\!-\!\lambda)\sqrt{2\theta})}\\
&= -\lim_{a\downarrow\lambda}
\, \frac{e^{(a-x)\sqrt{2\theta}}}{e^{(a-\lambda)\sqrt{2\theta}} }
\frac{\sinh((a-\lambda)\sqrt{2\theta})- \cosh((a-\lambda)\sqrt{2\theta})}{e^{(a-\lambda)\sqrt{2\theta}}}
 = e^{-(x-\lambda)\sqrt{2\theta}}
  = \E^x[e^{-\theta\tau_\lambda}].
\end{split}
\end{equation*}
In addition, by \cite[Sections II.6: displays 2.3.1, 2.0.2]{BorSal}, for $d\ge3$ it holds 
\begin{equation}\label{normalization_iden}
\begin{split}
&\;\lim_{\delta\downarrow0} \frac{1}{\delta^2} \int_\lambda^{\lambda+\delta} 
\E^x\big[L^a_{\tau_\lambda}\,\mathbf{1}_{\{\tau_\lambda<\infty\}}\big]\,\mathrm{d}a
 =\frac{1}{2}\lim_{a\downarrow\lambda}\,\frac{\mathrm{d}}{\mathrm{d}a}\, \E^x\big[L^a_{\tau_\lambda}\,\mathbf{1}_{\{\tau_\lambda<\infty\}} \big]\\
&= -\frac{1}{2}\lim_{a\downarrow\lambda}\,\lim_{b\downarrow0}
\, \frac{1}{b}\,\frac{\mathrm{d}}{\mathrm{d}a}\, \E^x\big[e^{- b L^a_{\tau_\lambda}}\,\mathbf{1}_{\{\tau_\lambda<\infty\}} \big]
= -\frac{1}{2}\lim_{a\downarrow\lambda}\,\lim_{b\downarrow0}
\, \frac{1}{b}\,\frac{\mathrm{d}}{\mathrm{d}a}\,
\frac{(d/2\!-\!1)\,x^{2-d}}{(d/2\!-\!1)\,\lambda^{2-d} \!+\! ba(\lambda^{2-d} \!-\! a^{2-d})}\\
&=\frac{x^{2-d}}{\lambda^{2-d}}=\PP^x(\tau_\lambda<\infty).
\end{split}
\end{equation}
For $d=1,2$, \cite[Sections II.3, II.4: displays 2.3.1]{BorSal} show that the first expression in \eqref{normalization_iden} is equal to $1$, and thus equal to the last expression in \eqref{normalization_iden}.

\medskip

In view of \eqref{lapl_trans_iden}, \eqref{normalization_iden}, and since the weak convergence of probability measures on $\RR_+$ follows from the pointwise convergence of their Laplace transforms (see~\cite[Example~5.5]{Bil}), 
\begin{equation*}
\frac{\gamma x^{d-1}}{\lambda}\,\lim_{\delta\downarrow0} \,\frac{1}{\delta^2}\, \int_{\lambda-\delta}^{\lambda+\delta} 
\E^x\big[L^a_{\tau_\lambda}\,\mathbf{1}_{\{\tau_\lambda\le s\}}\big]\,\mathrm{d}a = H(\lambda)\,\PP^x(\tau_\lambda\le s)\,x^{d-1}. 
\end{equation*}
(Note that $\PP^x(\tau_\lambda\!=\!s)\!=\!0$, $s\!\in\!\RR_+$, see \cite[Section II.6: display 2.02]{BorSal}.) The second summand in \eqref{loc_time_rep} vanishes since $\lim_{\delta\downarrow0} \frac{1}{\delta^2}\, \int_{\lambda-\delta}^{\lambda+\delta} 
L^a_s\,\mathbf{1}_{\{\tau_\lambda>s\}}\,\mathrm{d}a=0$ almost surely and
\begin{equation*}
\begin{split}
\E^x\bigg[\frac{1}{\delta^4}\bigg(\int_{\lambda-\delta}^{\lambda+\delta} 
L^a_s\,\mathbf{1}_{\{\tau_\lambda>s\}}\,\mathrm{d}a\bigg)^2\bigg] 
\le \E^x\bigg[\frac{1}{\delta^4}\bigg(\int_\lambda^{\lambda+\delta} 
L^a_{\tau_\lambda}\,\mathrm{d}a\bigg)^2\bigg] 
 \le \frac{1}{\delta^3}\,\E^x\bigg[\int_\lambda^{\lambda+\delta}  (L^a_{\tau_\lambda})^2\,\mathrm{d}a\bigg] \;\;\; \\
=\begin{cases}
\frac{1}{\delta^3}\int_\lambda^{\lambda+\delta} 8(a-\lambda)^2\,\mathrm{d}a\quad\text{if}\quad d=1,\\
\frac{1}{\delta^3}\int_\lambda^{\lambda+\delta} 8a^2(\log a-\log \lambda)^2\,\mathrm{d}a\quad\text{if}\quad d=2,\\
\frac{1}{\delta^3}\int_\lambda^{\lambda+\delta} 8(d-2)^{-2}\,x^{2-d}\,\lambda^{2d-4}\,a^d\,(a^{2-d}-\lambda^{2-d})^2\,\mathrm{d}a\quad\text{if}\quad d\ge3
\end{cases}
\end{split}
\end{equation*}
implies the uniform integrability of $\frac{1}{\delta^2}\, \int_{\lambda-\delta}^{\lambda+\delta} 
L^a_s\,\mathbf{1}_{\{\tau_\lambda>s\}}\,\mathrm{d}a$ as $\delta\downarrow0$. Here, we have used the monotonicity of local time, Jensen's inequality, and the first Ray-Knight Theorem (see, e.g., \cite[Chapter 6, Theorem 4.7]{KaSh}), \cite[Sections II.4, II.6: displays 2.3.2]{BorSal}, respectively.

\medskip

For $x<\lambda$, the first summand in \eqref{loc_time_rep} equals to $H(\lambda)\,\PP^x(\tau_\lambda\le s)\,x^{d-1}$ due to the analogues of \eqref{lapl_trans_iden}, \eqref{normalization_iden} based on \cite[Sections II.3, II.4, II.6:~displays 2.3.3, 2.0.1, 2.3.1, 2.0.2, and Appendix 2.10]{BorSal}.~The second summand in \eqref{loc_time_rep} is seen to be $0$ by a repetition of the same argument, where we now compute $\E^x[(L^a_{\tau_\lambda})^2]$ by the Ray-Knight Theorem in \cite[Theorem 2]{Yor}. All in all,
\begin{equation*}
\lim_{\delta\downarrow0}\,\sum_{i\ge1}\,\PP(Y^{\delta,i,\Delta}_t\in \cdot,\,(m-1)\Delta<T_i^{\delta,\Delta}\le t<\tau^{\delta,\Delta}_i)=H(\lambda)\,\PP^\cdot(\tau_\lambda\le s)\,\nu,
\end{equation*}
where the $\delta\downarrow0$ limit is taken in the sense of pointwise convergence of (Lebesgue) densities. 

\medskip

\noindent\textbf{Step 3.} Let us verify the applicability of the Dominated Convergence Theorem to the function on $\RR_+$ that maps $x$ to 
\begin{equation}\label{DCTintegrand}
\frac{\gamma x^{d-1}}{\lambda}\,\frac{1}{\delta^2}\, \int_{\lambda-\delta}^{\lambda+\delta} 
\E^x\big[L^a_{\tau_\lambda}\,\mathbf{1}_{\{\tau_\lambda\le s\}}\big]\,\mathrm{d}a 
+  \frac{\gamma x^{d-1}}{\lambda}\,\E^x\bigg[\frac{1}{\delta^2}\, \int_{\lambda-\delta}^{\lambda+\delta} 
L^a_s\,\mathbf{1}_{\{\tau_\lambda>s\}}\,\mathrm{d}a\bigg], 	
\end{equation}
as $\delta\downarrow0$.
To this end, we upper bound \eqref{DCTintegrand} by $\frac{x^{d-1}}{\lambda}\,\frac{\gamma}{\delta^2}\,\int_{\lambda-\delta}^{\lambda+\delta} \E^x\big[L^a_{\tau_\lambda}]\,\mathrm{d}a$ and use \cite[Theorem~2]{Yor} for $x<\lambda$, as well as \cite[Chapter 6, Theorem 4.7]{KaSh}, \cite[Sections II.4, II.6: displays 2.3.2]{BorSal} for $\lambda<x<\lambda+1$, to find that the latter equals to
\begin{equation*}
\begin{cases}
\frac{2\gamma}{\delta^2}\int_{\lambda-\delta}^\lambda (\lambda-a)\,\mathrm{d}a\quad\text{if}\quad d=1, \\
x\,\frac{2\gamma}{\delta^2}\int_{\lambda-\delta}^\lambda a(\log\lambda-\log a)\,\mathrm{d}a\quad\text{if}\quad d=2, \\
x^{d-1}\,\frac{2\gamma}{(d-2)\,\lambda^{d-2}\,\delta^2} \int_{\lambda-\delta}^\lambda a(\lambda^{d-2}-a^{d-2})\,\mathrm{d}a\quad\text{if}\quad d\ge 3
\end{cases}
\end{equation*}
and
\begin{equation*}
\begin{cases}
\frac{2\gamma}{\lambda\,\delta^2}\int_\lambda^{\lambda+\delta} (a-\lambda)\,\mathrm{d}a \quad\text{if}\quad d=1, \\
x\,\frac{2\gamma}{\lambda\,\delta^2}\int_\lambda^{\lambda+\delta} a(\log a-\log\lambda)\,\mathrm{d}a \quad\text{if}\quad d=2, \\
x\,\frac{2\gamma}{(d-2)\,\lambda^{5-2d}\,\delta^2} \int_\lambda^{\lambda+\delta} a(\lambda^{2-d}-a^{2-d})+a^{2d-3}(\lambda^{2-d}-a^{2-d})^2\,\mathrm{d}a\quad\text{if}\quad d\ge 3,
\end{cases}
\end{equation*}
respectively. For $x\ge\lambda+1$, we estimate \eqref{DCTintegrand} by $\frac{\gamma x^{d-1}}{\lambda}\,\PP^x(\tau_{\lambda+1}<s)\,\frac{1}{\delta^2}\,\int_\lambda^{\lambda+\delta} \E^{\lambda+1}\big[L^a_{\tau_\lambda}]\,\mathrm{d}a$. By \cite[Chapter 6, Theorem 4.7]{KaSh}, \cite[Sections II.4, II.6: displays 2.3.2]{BorSal}, this computes to
\begin{equation*}
\begin{cases}
\PP^x(\tau_{\lambda+1}<s)\,\frac{2\gamma}{\lambda\,\delta^2}\int_\lambda^{\lambda+\delta} (a-\lambda)\,\mathrm{d}a\quad\text{if}\quad d=1, \\
x\,\PP^x(\tau_{\lambda+1}<s)\,\frac{2\gamma}{\lambda\,\delta^2}\int_\lambda^{\lambda+\delta} a(\log a-\log\lambda)\,\mathrm{d}a\quad\text{if}\quad d=2, \\
x^{d-1}\,\PP^x(\tau_{\lambda+1}<s)\,\frac{2\gamma(\lambda+1)^{2-d}}{(d-2)\,\lambda^{5-2d}\,\delta^2}\int_\lambda^{\lambda+\delta}a(\lambda^{2-d}-a^{2-d})+a^{2d-3}(\lambda^{2-d}-a^{2-d})^2\,\mathrm{d}a\quad\text{if}\quad d\ge 3.
\end{cases}
\end{equation*}
Thus, the expression in \eqref{DCTintegrand} can be estimated by $C\,\mathbf{1}_{[0,\lambda+1)}(x)+Cx^{d-1}\PP^x(\tau_{\lambda+1}\!<\!s)\,\mathbf{1}_{[\lambda+1,\infty)}(x)$ for some $C=C(d,\gamma,\lambda)<\infty$.~A comparison with a standard Brownian motion started from~$x$ reveals the latter as (Lebesgue) integrable over $\RR_+$. Combining this conclusion with the results of Steps 1 and 2 we obtain \eqref{forward_to_backward} for $t\in((m-1)\Delta,m\Delta)$.

\medskip

\noindent\textbf{Step 4.} Finally, we note that the reasoning thus far also reveals $u^\Delta(m\Delta-,\cdot)\,\nu$ as
\begin{equation*}
\PP(X_{m\Delta}\!\in\!\cdot\,,\tau^\Delta\!\ge\!m\Delta)
\!+\!\!\sum_{n=1}^{m-1}\!
\sum_{i\ge1} \PP(X^{n,i,\Delta}_{m\Delta}\!\in\!\cdot\,,\tau^{n,\Delta}_i\!\ge\!m\Delta)
\!+\!\lim_{\delta\downarrow0} \sum_{i\ge1} \PP(Y^{\delta,i,\Delta}_{m\Delta}\!\in\!\cdot\,,T_i^{\delta,\Delta}\!<\! m\Delta\!\le\!\tau^{\delta,\Delta}_i),
\end{equation*}
where the $\delta\downarrow0$ limit is in the sense of the total variation convergence of measures. Together with the definition of $u^\Delta(m\Delta,\cdot)$ (Definition \ref{def:Euler}) this shows that $u^\Delta(m\Delta,\cdot)\,\nu$ equals to
\begin{equation}\label{uDeltaupdate}
\PP(X_{m\Delta}\!\in\!\cdot\,,\tau^\Delta\!>\!m\Delta)
\!+\!\!\sum_{n=1}^m
\sum_{i\ge1} \PP(X^{n,i,\Delta}_{m\Delta}\!\in\!\cdot\,,\tau^{n,\Delta}_i\!>\!m\Delta)
\!+\!\lim_{\delta\downarrow0} \sum_{i\ge1} \PP(Y^{\delta,i,\Delta}_{m\Delta}\!\in\!\cdot\,,T_i^{\delta,\Delta}\!\le\! m\Delta\!<\!\tau^{\delta,\Delta}_i),
\end{equation}
where the limit is taken in the same sense as before.~This completes the proof of \eqref{forward_to_backward}.

\medskip

The first equation in \eqref{Lambda_forward} follows from the induction hypothesis and
\begin{equation*}
\begin{split}
&\;\frac{1}{d}\big((\Lambda^\Delta_{m\Delta})^d-(\Lambda^\Delta_{(m-1)\Delta})^d\big)
= \int_{\RR_+} u^\Delta((m-1)\Delta,x)\,\nu(\mathrm{d}x) - \int_{\RR_+} u^\Delta(m\Delta,x)\,\nu(\mathrm{d}x)\\
&=-\PP(\tau^\Delta>m\Delta)+\PP(\tau^\Delta>(m-1)\Delta)
-\sum_{n=1}^m \sum_{i\ge1} \PP(\tau^{n,\Delta}_i> m\Delta) \\
&\; +\!\sum_{n=1}^{m-1}
\sum_{i\ge1} \PP(\tau^{n,\Delta}_i\!>\!(m\!-\!1)\Delta) 
\!-\!\lim_{\delta\downarrow0} \sum_{i\ge 1} \PP(T^{\delta,\Delta}_i\!\le\! m\Delta\!<\!\tau^{\delta,\Delta}_i)
\!+\!\lim_{\delta\downarrow0} \sum_{i\ge 1} 
\PP(T^{\delta,\Delta}_i\!\le\! (m\!-\!1)\Delta\!<\!\tau^{\delta,\Delta}_i),
\end{split}
\end{equation*}
which, in turn, is a consequence of Definition \ref{def:Euler}, \eqref{uDeltaupdate}, and the induction hypothesis. \qed

\begin{remark}\label{re:glob_bound}
The first equation in \eqref{Lambda_forward} implies $\sup_{\RR_+} \Lambda^\Delta\le(d+(\Lambda_{0-})^d)^{1/d}$.	
\end{remark}

\begin{lemma}\label{le:osc}
Fix a $\Delta>0$ and define 
\begin{align}
\sigma^\Delta=\inf\{t\in\Delta\NN\!: H(\Lambda^\Delta_t)\ge \|u(0-,\cdot)\|_\infty\}\wedge T.\label{eq.Sec2.sigmaDelta.def}
\end{align}
Then:
\begin{enumerate}[(a)]
\item $\Lambda^\Delta_{m_2\Delta}-\Lambda^\Delta_{m_1\Delta}\le C\sqrt{(m_2-m_1)\Delta}$ for all 
$0\le m_1\Delta<m_2\Delta<\sigma^\Delta$ with $(m_2-m_1)\Delta\le 1$, where $C=C(\Lambda_{0-},\|u(0-,\cdot)\|_\infty)<\infty$.
\item $\Lambda^\Delta$ is non-increasing on $[\sigma^\Delta,T]$.  	
\end{enumerate}
\end{lemma}

\noindent\textbf{Proof. (a).} Clearly, it suffices to consider $0\le m_1\Delta<m_2\Delta<\sigma^\Delta$ with $(m_2-m_1)\Delta\le 1$ and $\Lambda^\Delta_{m_1\Delta}=\inf_{[m_1\Delta,m_2\Delta]} \Lambda^\Delta=:\underline{\lambda}$, $\Lambda^\Delta_{m_2\Delta}=\sup_{[m_1\Delta,m_2\Delta]} \Lambda^\Delta=:\overline{\lambda}$. Then, the forward representation \eqref{Lambda_forward} gives
\begin{equation*}
\begin{split}
&\frac{1}{d}\big((\Lambda^\Delta_{m_2\Delta})^d\!-\!(\Lambda^\Delta_{m_1\Delta})^d\big) \!= \PP(\tau^\Delta\!\le\! m_2\Delta)
\!-\!\!\sum_{n=1}^{m_2} \sum_{i\ge 1} \PP(\tau^{n,\Delta}_i\!>\! m_2\Delta) 
\!-\! \lim_{\delta\downarrow0}\sum_{i\ge1} \PP(T^{\delta,\Delta}_i\!\le\! m_2\Delta\!<\!\tau^{\delta,\Delta}_i) \\
&\qquad\qquad\qquad\qquad\qquad\;\;
\!-\!\PP(\tau^\Delta\!\le\! m_1\Delta)
\!+\!\!\sum_{n=1}^{m_1} \sum_{i\ge 1} \PP(\tau^{n,\Delta}_i\!>\! m_1\Delta)
\!+\! \lim_{\delta\downarrow0}\sum_{i\ge1}
\PP(T^{\delta,\Delta}_i\!\le\! m_1\Delta\!<\!\tau^{\delta,\Delta}_i) \\
& \le\PP(m_1\Delta\!<\!\tau^\Delta\!\le\! m_2\Delta)
\!+\! \sum_{n=1}^{m_1} \sum_{i\ge1} \PP(m_1\Delta\!<\!\tau^{n,\Delta}_i\!\le \! m_2\Delta)
\!+\! \lim_{\delta\downarrow0}\,\sum_{i\ge1}\, \PP(T_i^{\delta,\Delta}
\!\le\! m_1\Delta \!<\! \tau^{\delta,\Delta}_i \!\le\! m_2\Delta) \\
&=\int_{\RR_+} \PP^x(\tau_{\Lambda^\Delta_{m_1\Delta+\cdot}}\le (m_2-m_1)\Delta)\, u^\Delta(m_1\Delta,x)\,\nu(\mathrm{d}x).
\end{split}
\end{equation*}
Recalling $\underline{\lambda}\!=\!\inf_{[m_1\Delta,m_2\Delta]} \Lambda^\Delta$, $\overline{\lambda}\!=\!\sup_{[m_1\Delta,m_2\Delta]} \Lambda^\Delta$ we bound the latter expression further~by
\begin{equation*}
\begin{split}
\int_{\underline{\lambda}}^{\overline{\lambda}} u^\Delta(m_1\Delta,x)\,\nu(\mathrm{d}x)
&+\int_{\overline{\lambda}}^\infty \PP^x(\tau_{\overline{\lambda}}\le(m_2-m_1)\Delta)\,u^\Delta(m_1\Delta,x)\,\nu(\mathrm{d}x) \\
&
+\int_0^{\underline{\lambda}} 
\PP^x(\tau_{\underline{\lambda}}\le(m_2-m_1)\Delta)\,u^\Delta(m_1\Delta,x)\,\nu(\mathrm{d}x).
\end{split}
\end{equation*}

\smallskip

Proceeding inductively over the intervals $[0,\Delta],\,[\Delta,2\Delta],\,\ldots$ we read off $u^\Delta\le\|u(0-,\cdot)\|_\infty$ on $[0,m_1\Delta]\times\RR_+$ from Definition \ref{def:Euler}. Plugging in this estimate we upper bound the above~by
\begin{equation*}
\|u(0-,\cdot)\|_\infty\bigg(\frac{1}{d}(\overline{\lambda}^d-\underline{\lambda}^d)
+\int_{\overline{\lambda}}^\infty \PP^x(\tau_{\overline{\lambda}}\le(m_2-m_1)\Delta)\,\nu(\mathrm{d}x)
+\int_0^{\underline{\lambda}} 
\PP^x(\tau_{\underline{\lambda}}\le(m_2-m_1)\Delta)\,\nu(\mathrm{d}x)\!\bigg).
\end{equation*}
Next, we bound the sum of the two integrals by replacing $\PP^x$ with the law of a standard Brownian motion started from $x$ when $x\ge\overline{\lambda}$, estimating $\PP^x(\tau_{\underline{\lambda}}\le(m_2-m_1)\Delta)$ by $1$ when $x\in[(\underline{\lambda}^2-d(m_2-m_1)\Delta)^{1/2},\underline{\lambda})=:[\lambda_*,\underline{\lambda})$, and applying the Dambis-Dubins-Schwarz Theorem (see, e.g., \cite[Chapter 3, Problem 4.7]{KaSh}) when $x\in[0,\lambda_*)$, upon noting that the diffusion coefficient of $R^2$ is smaller or equal to $4\underline{\lambda}^2$ until $\tau_{\underline{\lambda}}$ and that its drift coefficient is $d$.~Thus, the sum of the two integrals in the above display is at most
\begin{equation*}
\begin{split}
&\,\int_{\overline{\lambda}}^\infty 2\overline{\Phi}\bigg(\frac{x-\overline{\lambda}}{\sqrt{(m_2-m_1)\Delta}}\bigg)\,\nu(\mathrm{d}x)
+\int_{\lambda_*}^{\underline{\lambda}} \nu(\mathrm{d}x)
+\int_0^{\lambda_*} 2\overline{\Phi}\bigg(\frac{\underline{\lambda}^2-d(m_2-m_1)\Delta-x^2}{2\underline{\lambda}\sqrt{(m_2-m_1)\Delta}}\bigg)\,\nu(\mathrm{d}x) \\
& = 2\sqrt{(m_2-m_1)\Delta}\int_0^\infty \overline{\Phi}(y)(\overline{\lambda}+y\sqrt{(m_2-m_1)\Delta})^{d-1}\,\mathrm{d}y
+\frac{1}{d}\big(\underline{\lambda}^d-\lambda_*^d\big) \\
&\quad
+2\underline{\lambda}\sqrt{(m_2-m_1)\Delta}\int_0^{\frac{\lambda_*^2}{2\underline{\lambda}\sqrt{(m_2-m_1)\Delta}}} \overline{\Phi}(y)\,\big(\lambda_*^2-2y\underline{\lambda}\sqrt{(m_2-m_1)\Delta}\big)^{d/2-1}\,\mathrm{d}y,
\end{split}
\end{equation*}
where $\overline{\Phi}(x)$ is the probability that a standard normal random variable exceeds $x$.
In view of Remark \ref{re:glob_bound}, the above can be controlled by $C\sqrt{(m_2-m_1)\Delta}$ for some $C=C(\Lambda_{0-})<\infty$.

\medskip

All in all, we have obtained 
\begin{equation*}
\frac{1}{d}\big((\Lambda^\Delta_{m_2\Delta})^d-(\Lambda^\Delta_{m_1\Delta})^d\big)
\le \|u(0-,\cdot)\|_\infty\bigg(\frac{1}{d}(\overline{\lambda}^d-\underline{\lambda}^d)
+C\sqrt{(m_2-m_1)\Delta}\bigg),
\end{equation*}
which due to $\underline{\lambda}=\Lambda_{m_1\Delta}^\Delta$, $\overline{\lambda}=\Lambda_{m_2\Delta}^\Delta$ can be written as
\begin{equation*}
\frac{1}{d}\big((\Lambda^\Delta_{m_2\Delta})^d-(\Lambda^\Delta_{m_1\Delta})^d\big)
\le\frac{C\|u(0-,\cdot)\|_\infty}{1-\|u(0-,\cdot)\|_\infty}\,\sqrt{(m_2-m_1)\Delta}.
\end{equation*}
Since $\Lambda^\Delta_{m_1\Delta}\!>\!\gamma/\|u(0-,\cdot)\|_\infty$, the left-hand side bounds $(\Lambda^\Delta_{m_2\Delta}-\Lambda^\Delta_{m_1\Delta})\gamma^{d-1}/\|u(0-,\cdot)\|_\infty^{d-1}$.

\medskip 

\noindent\textbf{(b).} It suffices to check that, for $\sigma^\Delta<T$ and all $m=0,\,1,\,\ldots,\,\lfloor(T-\sigma^\Delta)/\Delta\rfloor$,
\begin{equation}\label{Lambda_decr_ind}
u^\Delta(\sigma^\Delta+m\Delta,\cdot)\le H(\Lambda^\Delta_{\sigma^\Delta+m\Delta})
\quad\text{and}\quad\Lambda^\Delta_{\sigma^\Delta+m\Delta}\le \Lambda^\Delta_{\sigma^\Delta+(m-1)^+\Delta}.  
\end{equation} 
To this end, we argue by induction over $m$. We have $H(\Lambda^\Delta_{\sigma^\Delta})\ge\|u(0-,\cdot)\|_\infty$ for $m\!=\!0$, whereas $H(\Lambda^\Delta_t)\!<\!\|u(0-,\cdot)\|_\infty$, $t\!\in\![0,\sigma^\Delta)$. Proceeding inductively over the intervals $[0,\Delta],\,[\Delta,2\Delta],\,\ldots$ we read off $u^\Delta\le\|u(0-,\cdot)\|_\infty$ on $[0,\sigma^\Delta)\times\RR_+$ from Definition \ref{def:Euler}. Hence, $u^\Delta(\sigma^\Delta,\cdot)\le H(\Lambda^\Delta_{\sigma^\Delta})$. Suppose now that \eqref{Lambda_decr_ind} holds for an $0\le m<\lfloor(T-\sigma^\Delta)/\Delta\rfloor$. Then, by Definition \ref{def:Euler}, the reversibility of $R$ with respect to $\nu$, and \eqref{Lambda_decr_ind},
\begin{equation*}
\begin{split}
& \,\int_0^\infty u^\Delta((\sigma^\Delta+(m+1)\Delta)-,x)\,\nu(\mathrm{d}x) \\
& =\int_0^\infty \E^x[u^\Delta(\sigma^\Delta\!+\!m\Delta,R_\Delta)]+ \E^x\big[\mathbf{1}_{\{\tau_{\Lambda^\Delta_{\sigma^\Delta+m\Delta}}<\Delta\}}\,
\big(H(\Lambda^\Delta_{\sigma^\Delta+m\Delta})
-u^\Delta(\sigma^\Delta\!+\!m\Delta,R_\Delta)\big)\big]\,\nu(\mathrm{d}x) \\
&\ge \int_0^\infty \E^x[u^\Delta(\sigma^\Delta\!+\!m\Delta,R_\Delta)]\,\nu(\mathrm{d}x)
=\int_0^\infty u^\Delta(\sigma^\Delta+m\Delta,x)\,\nu(\mathrm{d}x).
\end{split}
\end{equation*}
Therefore, Definition \ref{def:Euler} yields $\Lambda^\Delta_{\sigma^\Delta+(m+1)\Delta}\le \Lambda^\Delta_{\sigma^\Delta+m\Delta}$. Using Definition \ref{def:Euler}, \eqref{Lambda_decr_ind}, and the latter inequality we find
\begin{equation*}
u^\Delta(\sigma^\Delta+(m+1)\Delta,\cdot)
\le u^\Delta(\sigma^\Delta+(m+1)\Delta-,\cdot)
\vee H(\Lambda^\Delta_{\sigma^\Delta+(m+1)\Delta})
\le H(\Lambda^\Delta_{\sigma^\Delta+(m+1)\Delta}),
\end{equation*}
thus concluding the proof. \qed

\medskip

We are now ready to give the proof of Proposition \ref{prop:tightness}.

\medskip

\noindent\textbf{Proof of Proposition \ref{prop:tightness}.} Remark \ref{re:glob_bound} and Lemma \ref{le:osc} allow to easily verify the relative compactness criterion for $D([-1,T+1],\RR)$ with the M1 topology (see, e.g., \cite[Theorem 4.3]{DIRT2}). \qed

\medskip

\begin{remark}\label{rem:uDelta.bound}
Employing the arguments in the proof of Lemma \ref{le:osc} it is easy to deduce from Assumption \ref{main_ass} and Definition \ref{def:Euler} the existence of a constant $C_2<\infty$ such that
\begin{align*}
&u^\Delta(t,x)\leq C_1(t)\wedge \left(C_2(C_1(t)+1) e^{-\alpha x/2}\right),\quad x>0,\quad\text{with} \\
& C_1(t):= \|u(0-,\cdot)\|_\infty\vee \sup_{s\in[0,t]}H(\Lambda^\Delta_{s}),
\end{align*}
holds for all $t\in[0,T]$. This observation is used in Subsection \ref{subse:upperBound}.
\end{remark}

\section{Limit points as solutions of the Stefan problem}
\label{se:3}

In this section, we prove that every limit point $\Lambda$ of $\{\Lambda^\Delta\}_{\Delta>0}$ solves the Stefan problem. For any $\Lambda\in D([-1,T],\RR_+)$, with $\Lambda_{0-}>0$, we recall the associated $\zeta=\inf\{t>0\!: \Lambda_t=0\}$ and consider
\begin{equation}\label{what_is_u}
u(t,x)=\E^x\big[\mathbf{1}_{\{\tau_{\Lambda_{t-\cdot}}\leq t\}}\,H(R_{\tau_{\Lambda_{t-\cdot}}})\big]
+\E^x\big[\mathbf{1}_{\{\tau_{\Lambda_{t-\cdot}}> t\}}\,u(0-,R_t)\big],\quad (t,x)\in[0,T\wedge\zeta]\times(0,\infty),
\end{equation}
where $\tau_{\Lambda_{t-\cdot}}=\inf\{s\ge0\!:(R_s-\Lambda_{t-s})(x-\Lambda_t)\le 0\}$ as before.
For convenience, we set $u(t,x):=u(0-,x)$ for $t<0$.

\medskip

Let us establish several preliminary properties of $u$. First, the Markov property of the Bessel process yields that, under Assumption \ref{main_ass}, any $u$ satisfying \eqref{what_is_u} also satisfies for any $0\leq s<t\leq T\wedge\zeta$,
\begin{equation}\label{eq.Sec4.FK.t.s}
u(t,x)=\E^x\big[\mathbf{1}_{\{\tau_{\Lambda_{t-\cdot}}\leq t-s\}}\,H(R_{\tau_{\Lambda_{t-\cdot}}})\big]+\E^x\big[\mathbf{1}_{\{\tau_{\Lambda_{t-\cdot}}> t-s\}}\,u(s,R_{t-s})\big], \quad x>0.
\end{equation}



\smallskip

Next, we establish an upper bound on $u$.

\begin{lemma}\label{le:Sec4.u.bound.H}
Let $u$ be given by \eqref{what_is_u}. 
Then, for any $t\in[0,T\wedge\zeta]$, we have 
\begin{align*}
&u(t,x)\leq C_1(t,x)\wedge \left(C_2(C_1(t,x)+1) e^{-\alpha x/2}\right),\quad x>0,\quad\text{with} \\
& C_1(t,x):= \|u(0-,\cdot)\|_\infty\vee \sup_{s\in[0,t]}H(\Lambda_{s-})\vee(H(x)\wedge H(\Lambda_t)),
\end{align*}
where $C_2<\infty$ is a constant and $H(0):=\infty$.
\end{lemma}
\noindent\textbf{Proof.}
The upper bound $\|u(0-,\cdot)\|_\infty\vee \sup_{s\in[0,t]}H(\Lambda_{s-})\vee \left(H(x)\wedge H(\Lambda_t)\right)$ follows directly from \eqref{what_is_u}.~To obtain the upper bound $C_2(C_1(t,x)+1) e^{-\alpha x/2}$, we combine \eqref{what_is_u}, the fact that the $d$-dimensional Bessel process is the radial part of a Brownian motion in $\RR^d$ and \cite[displays (18), (19)]{Law}, which yield that (i) the first term on the right-hand side of \eqref{what_is_u} decays faster than exponentially as $x\rightarrow\infty$ and that (ii) the second term can be bounded by $C e^{-\alpha x/2}$ (recall that $u(0-,x)\leq C e^{-\alpha x}$).
\qed

\medskip

The next lemma characterizes $u(t-,\cdot)$ and $u(t+,\cdot)$.

\begin{lemma}\label{le:u.LRlims}
Let $u$ be given by \eqref{what_is_u}.
Then, for $x\in\RR_+\setminus \{\Lambda_{t-},\Lambda_t\}$, $\lim_{s\uparrow t} u(s,x)=:u(t-,x)$ is well-defined for any $t\in[0,T\wedge\zeta]$ with $\Lambda_{t-}>0$, and $\lim_{s\downarrow t} u(s,x)=:u(t+,x)$ is well-defined for any $t\in[0,T\wedge\zeta)$.
In addition, 
\begin{align}
& u(t+,x)=u(t,x),\quad x\in\RR_+\setminus \{\Lambda_{t-},\Lambda_t\},\label{eq.Sec4.RightLim.u.goal}\\
& u(t-,x)=u(t,x),\quad x\in\RR_+\setminus[\Lambda_t\wedge\Lambda_{t-},\Lambda_t\vee\Lambda_{t-}],\label{eq.Sec4.LeftLim.u.goal}\\
& u(t,x)=H(x),\quad x\in(\Lambda_t\wedge\Lambda_{t-},\Lambda_t\vee\Lambda_{t-}).\label{eq.Sec4.inJump.u.goal}
\end{align}
\end{lemma}
\noindent\textbf{Proof.}
Let us prove that, for $x\notin \{\Lambda_{t-},\Lambda_t\}$ and as $s\uparrow t$, $u(s,x)$ converges to
\begin{equation}\label{eq.Sec4.LeftLim.u.1}
\E^x\big[\mathbf{1}_{\{\tau^{t-}\leq t\}}\,H(R_{\tau^{t-}})+\mathbf{1}_{\{\tau^{t-}> t\}}\,u(0-,R_t)\big],
\end{equation}
where $\tau^{t-}:=\inf\{s>0\!:(R_s-\Lambda_{t-s})(x-\Lambda_{t-})\le 0\}$.
To this end, we recall from \eqref{what_is_u} that 
\begin{equation}\label{eq.Sec4.LeftLim.u.2}
u(s,x)=\E^x\big[\mathbf{1}_{\{\tau_{\Lambda_{s-\cdot}}\leq s\}}\,H(R_{\tau_{\Lambda_{s-\cdot}}})
+\mathbf{1}_{\{\tau_{\Lambda_{s-\cdot}}> s\}}\,u(0-,R_s)\big],
\end{equation}
where $\tau_{\Lambda_{s-\cdot}}=\inf\{q\ge0\!:(R_q-\Lambda_{s-q})(x-\Lambda_s)\le 0\}$ as before.

\medskip

Let us show that the random variable inside the expectation in \eqref{eq.Sec4.LeftLim.u.2} converges almost surely to the one inside the expectation in \eqref{eq.Sec4.LeftLim.u.1} as $s\uparrow t$.
On the event $\{\tau^{t-}> t\}$, it is easy to see from the continuity of $R$ and the right-continuity of $\Lambda$ that, almost surely, $\tau_{\Lambda_{s-\cdot}}>t$ for all large enough $s\!\in\!(0,t)$. Similarly, on the event $\{\tau^{t-}\!=\! t\}$, we must have $R_t\!\in\!(\Lambda_0\wedge\Lambda_{0-},\Lambda_0\vee\Lambda_{0-})$ and $\inf_{q\in[0,t)} |R_q-\Lambda_{t-q}|>0$, from which we deduce that $R_s\in(\Lambda_0\wedge\Lambda_{0-},\Lambda_0\vee\Lambda_{0-})$ and $\tau_{\Lambda_{s-\cdot}}\!=\!s$ for all large enough $s\in(0,t)$.
Next, consider the event $\{\tau^{t-}< t\}$. First, we notice that $\tau^{t-}>0$ since $x\notin \{\Lambda_{t-},\Lambda_t\}$.
Then, we use the strong Markov property of $R$ and \cite[proof of Proposition 3.4(i)]{ABPR} to deduce that $R_\cdot - \Lambda_{t-\cdot}$ cannot have zeroes in $(0,t)$ that are local extrema. The latter conclusion and the left-continuity of $R_\cdot - \Lambda_{t-\cdot}$ imply that any small enough open right neighborhood of $\tau^{t-}$ contains a point $q$ such that $(R_q-\Lambda_{t-q})(x-\Lambda_{t-})< 0$ and $\Lambda_{t-q}=\Lambda_{(t-q)-}$. Then, for all large enough $s\in(0,t)$, we have $(R_q-\Lambda_{s-q})(x-\Lambda_{s})< 0$, which implies that $\overline{\lim}_{s\uparrow t}\tau_{\Lambda_{s-\cdot}}\leq \tau^{t-}$. On the other hand, using the continuity of $R$ and the right-continuity of $\Lambda$, we deduce, as above, that $\underline{\lim}_{s\downarrow t}\tau_{\Lambda_{s-\cdot}}\geq \tau^{t-}$, which altogether yields $\lim_{s\downarrow t}\tau_{\Lambda_{s-\cdot}}= \tau^{t-}<t$. Collecting these observations and employing the continuity of $H$ and $u(0-,\cdot)$, we conclude that the random variable inside the expectation in \eqref{eq.Sec4.LeftLim.u.2} converges almost surely to the one inside the expectation in \eqref{eq.Sec4.LeftLim.u.1} as $s\uparrow t$.

\medskip

It only remains to notice that $x\wedge\inf_{q\in[0,s]}\Lambda_{q-}\leq R_{\tau_{\Lambda_{s-\cdot}}}$ and to apply the Dominated Convergence Theorem (recall Assumption \ref{main_ass}), to deduce that the right-hand side of \eqref{eq.Sec4.LeftLim.u.2} converges to the expression in \eqref{eq.Sec4.LeftLim.u.1}. To obtain \eqref{eq.Sec4.LeftLim.u.goal}, we simply notice that $\tau^{t-}=\tau_{\Lambda_{t-\cdot}}$ for $x\in\RR_+\setminus[\Lambda_t\wedge\Lambda_{t-},\Lambda_t\vee\Lambda_{t-}]$. The convergence of $u(s,x)$ as $s\downarrow t$, and \eqref{eq.Sec4.RightLim.u.goal} follow from very similar arguments.
Finally, we deduce \eqref{eq.Sec4.inJump.u.goal} from \eqref{what_is_u} by observing that $\tau_{\Lambda_{t-\cdot}}=0$ for $x\in(\Lambda_t\wedge\Lambda_{t-},\Lambda_t\vee\Lambda_{t-})$.
\qed

\medskip

The next proposition is then the main result of this section.

\begin{proposition}\label{prop:Stefan}
Let $\Lambda$ be a limit point of $\{\Lambda^\Delta\}_{\Delta>0}$ in the sense of Proposition \ref{prop:tightness}, and let $u$ be defined by \eqref{what_is_u}. Then, the weak Stefan growth condition holds: 
\begin{equation}\label{growth}
\frac{1}{d}\big((\Lambda_t)^d-(\Lambda_{0-})^d\big)
=\int_{\RR_+} u(0-,x)\,\nu(\mathrm{d}x)-\int_{\RR_+} u(t,x)\,\nu(\mathrm{d}x),\quad t\in[0,T\wedge\zeta).
\end{equation} 
\end{proposition}

Our proof of Proposition \ref{prop:Stefan} is based on the following key lemma.

\begin{lemma}\label{lem:u_conv}
Let $\Delta_k\downarrow0$ be such that $\Lambda=\lim_{k\to\infty} \Lambda^{\Delta_k}$ in the sense of Proposition \ref{prop:tightness}. Then, with  $u$ defined by \eqref{what_is_u}, with $u^{\Delta_k}$ as in Definition \ref{def:Euler}, and with $\sigma_1=\sigma_1(t,k):=\sup\{s\le t\!: s\in\Delta_k\NN\}$, we have $u(t,x)=\lim_{k\to\infty} u^{\Delta_k}(\sigma_1,x)$ for all $t\in[0,T\wedge\zeta)$ satisfying $\Lambda_{t-}=\Lambda_t$ and all $\RR_+\ni x\neq\Lambda_t$.
\end{lemma}

\noindent\textbf{Proof.} Consider a $t\in[0,T\wedge\zeta)$ as described and an $x>\Lambda_t$. (The proof for $\RR_+\ni x<\Lambda_t$ is very similar.) We study $u(t,x)-u^{\Delta_k}(\sigma_1,x)$ via a coupling, i.e., by expressing it as a single expectation.~To this end, we let $R^x$ be a $d$-dimensional Bessel process started from $x$ and set
\begin{equation*}
\tau_k=\inf\{s\ge0\!:R^x_s\le \Lambda^{\Delta_k}_{\sigma_1-s}\}\wedge(\sigma_1+1)
\quad\text{and}\quad
\tau=\inf\{s\ge0\!:R^x_s\le \Lambda_{t-s}\}\wedge(t+1).
\end{equation*}
Then, using \eqref{what_is_u}, Definition \ref{def:Euler}, the Markov property of $R^x$, and $\Lambda^{\Delta_k}_{0-}=\Lambda^{\Delta_k}_0$ we infer 
\begin{equation}\label{eq:coupling}
u(t,x)-u^{\Delta_k}(\sigma_1,x)
=\E\big[\mathbf{1}_{\{\tau\leq t\}}H(R^x_\tau)
+\mathbf{1}_{\{\tau> t\}} u(0-,R^x_t)
-\mathbf{1}_{\{\tau_k\leq\sigma_1\}} H(R^x_{\tau_k})
-\mathbf{1}_{\{\tau_k>\sigma_1\}} u(0-,R^x_{\sigma_1})\big].
\end{equation}
On the event $\{\tau> t\}$, we have $R^x_t>\Lambda_0\wedge\Lambda_{0-}$ almost surely, so $\sup_{s\in[0,t+\varepsilon]} (\Lambda_{t-s}-R^x_s)<0$, almost surely, for a sufficiently small $\varepsilon\in(0,1)$.~Moreover, $s\mapsto\Lambda_s-R^x_{t-s}$ is the M1 limit of $s\mapsto\Lambda^{\Delta_k}_{\sigma_1-t+s}-R^x_{t-s}$ on $[-1,t]$, by definition of the M1 distance (e.g., \cite[Subsection 4.1]{DIRT2}) and \cite[Corollary 12.7.1]{Whi}. Thus, \cite[display (5.31)]{DIRT2} yields $\overline{\lim}_{k\to\infty} \sup_{s\in[-\varepsilon,t]} (\Lambda^{\Delta_k}_{\sigma_1-s}-R^x_s) \!\le\! \sup_{s\in[-\varepsilon,t]} (\Lambda_{t-s}-R^x_s)\!<\!0$.~It follows that $\tau_k> t\ge\sigma_1$ for all large enough $k$, almost surely. So, the random variable inside the expectation in \eqref{eq:coupling} is $u(0-,R^x_t)-u(0-,R^x_{\sigma_1})$ for all large enough $k$, almost surely.

\medskip

On the event $\{\tau=t\}$, we have, almost surely, $\Lambda_0-R^x_t<0<\Lambda_{0-}-R^x_t$, as well as $\sup_{s\in[0,t]}(\Lambda_s-R^x_{t-s})<0$. Then, the M1 convergence of $s\mapsto\Lambda^{\Delta_k}_{\sigma_1-t+s}-R^x_{t-s}$ to $s\mapsto\Lambda_s-R^x_{t-s}$, and $\Lambda_s=\Lambda^{\Delta_k}_s$ for $s\in[-1,0)$, imply that $\lim_{k\rightarrow\infty}\tau_k=\tau=t$ and that $\tau_k\leq \sigma_1$ for all large enough $k$.~The random variable inside the expectation in \eqref{eq:coupling} then equals to $H(R^x_\tau)-H(R^x_{\tau_k})$ almost surely, where $\lim_{k\to\infty} \tau_k=\tau$ almost surely. 

\medskip

We turn to the event $\{\tau<t\}$.~Repeating the argument in the preceding paragraphs we find $\underline{\lim}_{k\to\infty} \tau_k\ge\tau-\upsilon$ almost surely, for all $\upsilon>0$, and thus $\underline{\lim}_{k\to\infty} \tau_k\ge\tau$ almost surely. Further, the strong Markov property of $R^x$, Girsanov's Theorem, \cite[proof of Proposition 3.4(i)]{ABPR}, and \cite[Corollary 12.2.1]{Whi} almost surely allow us to find $t-\tau>\upsilon_l\downarrow 0$ such that $R^x_{\tau+\upsilon_l}<\Lambda_{t-\tau-\upsilon_l}$ and $\Lambda_{t-\tau-\upsilon_l} = \Lambda_{(t-\tau-\upsilon_l)-}$. Since $\Lambda_{t-\tau-\upsilon_l}=\lim_{k\to\infty} \Lambda^{\Delta_k}_{\sigma_1-\tau-\upsilon_l}$, we conclude that $R^x_{\tau+\upsilon_l}<\Lambda^{\Delta_k}_{\sigma_1-\tau-\upsilon_l}$, and consequently $\tau_k\le\tau+\upsilon_l<\sigma_1$, for all large enough~$k$. As above, we conclude that the random variable inside the expectation in \eqref{eq:coupling} equals to $H(R^x_\tau)-H(R^x_{\tau_k})$ almost surely, where $\lim_{k\to\infty} \tau_k=\tau$ almost surely.  

\smallskip

At this point, we see from the almost sure continuity of $R^x$ and the continuity of $u(0-,\cdot)$ and $H$ that the random variable  inside the expectation in \eqref{eq:coupling}  tends to $0$ almost surely in the limit $k\to\infty$. By applying the Dominated Convergence Theorem with the upper bound $H(\Lambda_{0-}\wedge\inf_{[0,t]}\Lambda)\vee\|u(0-,\cdot)\|_\infty\vee\sup_k H(\inf_{[0,t]} \Lambda^{\Delta_k})$ we get $u(t,x)=\lim_{k\to\infty} u^{\Delta_k}(\sigma_1,x)$. \qed

\medskip

We are now ready to give the proof of Proposition \ref{prop:Stefan}. 

\medskip

\noindent\textbf{Proof of Proposition \ref{prop:Stefan}.}\!\! Let $t\in[0,T\wedge\zeta)$ be such that $\Lambda_{t-}=\Lambda_t$.~With $\sigma_1$ as in Lemma~\ref{lem:u_conv}, our starting point is the identity
\begin{equation}\label{prel_cons}
\frac{1}{d}\big((\Lambda^{\Delta_k}_{\sigma_1})^d-(\Lambda_{0-})^d\big)
=\int_{\RR_+} u(0-,x)\,\nu(\mathrm{d}x)-\int_{\RR_+} u^{\Delta_k}(\sigma_1,x)\,\nu(\mathrm{d}x),
\end{equation}
which follows directly from Definition \ref{def:Euler}.~We take the limit $k\to\infty$ on both sides.~The left-hand side converges to $\frac{1}{d}\big((\Lambda_t)^d-(\Lambda_{0-})^d\big)$ thanks to the M1 convergence of $\Lambda^{\Delta_k}$ to $\Lambda$. On the right-hand side, Lemma \ref{lem:u_conv} gives $\lim_{k\to\infty} u^{\Delta_k}(\sigma_1,x)=u(t,x)$. Moreover, as in \eqref{eq:coupling},  
\begin{equation*}
\begin{split}
u^{\Delta_k}(\sigma_1,x)&=\E\big[\mathbf{1}_{\{\tau_k\leq\sigma_1\}}\,H(R^x_{\tau_k})
+\mathbf{1}_{\{\tau_k>\sigma_1\}}\,u(0-,R^x_{\sigma_1})\big] \\
&\le\sup_k H(\inf_{[0,t]} \Lambda^{\Delta_k})\,\PP(\tau_k\leq\sigma_1)+\E[u(0-,R^x_{\sigma_1})].
\end{split}
\end{equation*}
Estimating $\PP(\tau_k\leq\sigma_1)$ via Remark \ref{re:glob_bound} and a comparison of $R^x$ with a standard Brownian motion started from $x$, and controlling $\E[u(0-,R^x_{\sigma_1})]$ via $u(0-,y)\le e^{-\alpha y}$, $y\in\RR_+$ and \cite[display (18)]{Law}, we justify the application of the Dominated Convergence Theorem to the right-hand side of \eqref{prel_cons}. All in all, the $k\to\infty$ limit of \eqref{prel_cons} reads
\begin{equation*}
\frac{1}{d}\big((\Lambda_t)^d-(\Lambda_{0-})^d\big)
=\int_{\RR_+} u(0-,x)\,\nu(\mathrm{d}x)-\int_{\RR_+} u(t,x)\,\nu(\mathrm{d}x).
\end{equation*}

\smallskip

Next, let $t\in[0,T\wedge\zeta)$ be arbitrary. Then, there exist $T\wedge\zeta>t_l\downarrow t$ with $\Lambda_{t_l-}=\Lambda_{t_l}$ (see \cite[Corollary 12.2.1]{Whi}). In particular, 
\begin{equation}\label{growth_cont_points}
\frac{1}{d}\big((\Lambda_{t_l})^d-(\Lambda_{0-})^d\big)
=\int_{\RR_+} u(0-,x)\,\nu(\mathrm{d}x)-\int_{\RR_+} u(t_l,x)\,\nu(\mathrm{d}x).
\end{equation}
We take the limit $l\to\infty$ on both sides.~The left-hand side converges to $\frac{1}{d}\big((\Lambda_t)^d-(\Lambda_{0-})^d\big)$ thanks to the right-continuity of $\Lambda$.~On the right-hand side, we 
use Lemma \ref{le:u.LRlims} to deduce $\lim_{l\to\infty} u(t_l,x)=u(t,x)$ for all $x\in\RR_+\setminus\{\Lambda_{t-},\Lambda_t\}$.
Recalling the $\nu$-integrable upper bound on $u(t_l,\cdot)$ from Lemma \ref{le:Sec4.u.bound.H}, we pass to the $l\to\infty$ limit of \eqref{growth_cont_points} to find \eqref{growth}. \qed

\section{Limit points as physical solutions of the Stefan problem}
\label{se:4}

The goal of this section is to show that (i) any solution $(\Lambda,u)$ of \eqref{what_is_u}, \eqref{growth} with $d\geq3$ satisfies 
\begin{equation}\label{jump_cond}
\Lambda_{t-}-\Lambda_t \geq \inf\bigg\{y\in(0,\Lambda_{t-}):\, \int_{\Lambda_{t-}-y}^{\Lambda_{t-}} u(t-,x)\,\nu(\mathrm{d}x) >\int_{\Lambda_{t-}-y}^{\Lambda_{t-}} \big(H(x)-1\big)\,\nu(\mathrm{d}x)\bigg\},
\end{equation} 
and that (ii) any limit point $\Lambda$ of $\{\Lambda^\Delta\}_{\Delta\downarrow0}$, with $u$ defined via \eqref{what_is_u}, satisfies the above inequality with an equality.~This shows that the limits points $\Lambda$ have the smallest jumps among all solutions to \eqref{what_is_u}, \eqref{growth} and explains why we call such solutions \textit{physical}.

\subsection{Preliminary results}
\label{subse:Sec4.prelim}

The following lemma asserts that $\Lambda$ cannot jump at time $t$ if $\|u(t-,\cdot)\|_\infty\leq1$ and $\sup_{s\in[-1,t]}H(\Lambda_s)\leq1$. 

\begin{lemma}\label{le:noJumps.smallU}
Let $(\Lambda,u)$ be a solution to \eqref{what_is_u}, \eqref{growth}, and let $t\in[0,T\wedge\zeta)$ be such that $\|u(t-,\cdot)\|_\infty\leq 1$ and $\sup_{s\in[-1,t]}H(\Lambda_s)\leq1$. Then, $\Lambda_{t-}=\Lambda_t$.
\end{lemma}

\noindent\textbf{Proof.}
The growth condition \eqref{growth} yields, for $-1\leq s<t$:
\begin{equation*}
\frac{1}{d}\big((\Lambda_t)^d-(\Lambda_s)^d\big)
=\int_{\RR_+} u(s,x)-u(t,x)\,\nu(\mathrm{d}x).
\end{equation*}
Using the Dominated Convergence Theorem (recall Lemma \ref{le:Sec4.u.bound.H}) we pass to the limit in the above as $s\uparrow t$, and get
\begin{align*}
&\frac{1}{d}\big|(\Lambda_t)^d-(\Lambda_{t-})^d\big|
=\bigg|\int_{\RR_+} u(t-,x)-u(t,x)\,\nu(\mathrm{d}x)\bigg|
=\bigg|\int_{\Lambda_t\wedge\Lambda_{t-}}^{\Lambda_t\vee\Lambda_{t-}} u(t-,x)-H(x)\,\nu(\mathrm{d}x)\bigg|,
\end{align*}
where we also used Lemma \ref{le:u.LRlims} to obtain the last equality. It is easy to see that, if $\Lambda_{t-}\neq\Lambda_t$, the right-hand side in the above is strictly less than the left-hand side, which yields the conclusion of the lemma.
\qed

\medskip

An immediate corollary of Lemmas \ref{le:Sec4.u.bound.H} and \ref{le:noJumps.smallU} is that $\Lambda$ cannot have jumps before 
\begin{align}
\sigma:=\inf\big\{t>0:\,H(\Lambda_t)> \|u(0-,\cdot)\|_\infty\big\}.\label{eq.Sec4.sigma.def}
\end{align}
It is also clear that, if $\Lambda$ jumps at $\sigma$, the jump must be downwards.~The next lemma shows that $\Lambda$ is monotone after $\sigma$.

\begin{lemma}
Let $(\Lambda,u)$ be a solution to \eqref{what_is_u}, \eqref{growth} whose upward jumps satisfy 
\begin{align}\label{eq.Sec4.DeltaLambda.upper}
\Lambda_{t}-\Lambda_{t-} \leq \inf\bigg\{y>0:\, \int_{\Lambda_{t-}}^{\Lambda_{t-}+y} u(t-,x)\,\nu(\mathrm{d}x) <\int_{\Lambda_{t-}}^{\Lambda_{t-}+y} \big(H(x)+1\big)\,\nu(\mathrm{d}x)\bigg\}.
\end{align}
Then, $\Lambda$ is non-increasing on $[\sigma,T\wedge\zeta)$.
\end{lemma}

\noindent\textbf{Proof.} For $\eta\geq \|u(0-,\cdot)\|_\infty$, define
\begin{align*}
&\widetilde{\sigma}^\eta=\inf\big\{t>0:\,\sup_{s\in[0,t]}H(\Lambda_s)\geq \eta\big\},
\quad \sigma^\eta=\inf\{t\geq\widetilde{\sigma}^\eta:\,H(\Lambda_t)\neq H(\Lambda_{\widetilde{\sigma}^\eta})\}.
\end{align*}
Notice that Lemma \ref{le:Sec4.u.bound.H} and \eqref{eq.Sec4.DeltaLambda.upper} imply that, whenever $\sigma^\eta<T\wedge\zeta$, we have $\Lambda_{\sigma^\eta}\leq \Lambda_{\sigma^\eta-}$.
To prove the lemma, we argue by contradiction.~Assume that $\Lambda$ is not non-increasing on $[\sigma,T\wedge\zeta)$.~Then, there exists $\eta\geq H(\Lambda_{\sigma})\geq \|u(0-,\cdot)\|_\infty$ such that $\sigma^\eta<T\wedge\zeta$, $\Lambda_{\sigma^\eta}= \Lambda_{\sigma^\eta-}$, $\|u(\sigma^\eta,\cdot)\|_\infty \leq H(\Lambda_{\sigma^\eta})$, and for any $\varepsilon\in(0,T\wedge\zeta-\sigma^\eta)$ there exists an $\varepsilon'\in(0,\varepsilon)$ such that
\begin{align*}
& \Lambda_{\sigma^\eta} < \Lambda_{\sigma^\eta+\varepsilon'},
\quad\Lambda_{\sigma^\eta}\leq \Lambda_s\leq \Lambda_{\sigma^\eta+\varepsilon'},
\;\; s\in[\sigma^\eta,\sigma^\eta+\varepsilon'],
\quad\sup_{s\in[\sigma^\eta,\sigma^\eta+\varepsilon']}\|u(s,\cdot)\|_\infty=: C_1 <\infty.
\end{align*}
Next, using \eqref{growth}, \eqref{eq.Sec4.FK.t.s}, and setting $\tau=\tau_{\Lambda_{(\sigma^\eta+\varepsilon')-\cdot}}$ we obtain
\begin{align*}
&\;0<\frac{1}{d}\left((\Lambda_{\sigma^\eta+\varepsilon'})^d - (\Lambda_{\sigma^\eta})^d \right)
= \int_{\RR_+} u(\sigma^\eta,x)\,\nu(\mathrm{d}x) - \int_{\RR_+} u(\sigma^\eta+\varepsilon',x)\,\nu(\mathrm{d}x) \\
&= \int_{\RR_+} \EE^x\left[\bone_{\{\tau\leq\varepsilon'\}}\left(u(\sigma^\eta,R_{\varepsilon'}) - H(R_{\tau})\right) \right]\,\nu(\mathrm{d}x)
\leq \int_{\RR_+} \EE^x\left[\bone_{\{\tau\leq\varepsilon'\}}\left(H(\Lambda_{\sigma^\eta}) - H(R_{\tau})\right) \right]\,\nu(\mathrm{d}x) \\
&\leq C_2(\Lambda_{\sigma^\eta+\varepsilon'}-\Lambda_{\sigma^\eta}) \int_{\RR_+} \PP^x\left(\tau\leq\varepsilon'\right)\,\nu(\mathrm{d}x)\\
&\leq C_3(\Lambda_{\sigma^\eta+\varepsilon'}-\Lambda_{\sigma^\eta}) \bigg[ (\Lambda_{\sigma^\eta+\varepsilon'})^d-(\Lambda_{\sigma^\eta})^d
+ \int_{0}^{\Lambda_{\sigma^\eta}} \PP^x(\tau\leq\varepsilon')\,\nu(\mathrm{d}x) 
+ \int_{\Lambda_{\sigma^\eta+\varepsilon'}}^\infty \PP^x(\tau\leq\varepsilon')\,\nu(\mathrm{d}x)\bigg] \\
&\leq C_3(\Lambda_{\sigma^\eta+\varepsilon'}-\Lambda_{\sigma^\eta}) \bigg[ (\Lambda_{\sigma^\eta+\varepsilon'})^d-(\Lambda_{\sigma^\eta})^d
+ 2d\int_{0}^{\Lambda_{\sigma^\eta}} \PP\Big(\sup_{s\in[0,\varepsilon']} B_s\geq \frac{\Lambda_{\sigma^\eta}-x}{\sqrt{d}}\Big)\,\nu(\mathrm{d}x) \\ 
&\qquad\qquad\qquad\qquad\quad + \int_{\RR_+} \PP\Big(\inf_{s\in[0,\varepsilon']} B_s\leq -x\Big)(x+\Lambda_{\sigma^\eta+\varepsilon'})^{d-1}\,\mathrm{d}x\bigg],
\end{align*}
where $B$ is a standard Brownian motion.~It remains to divide the above by $\Lambda_{\sigma^\eta+\varepsilon'}-\Lambda_{\sigma^\eta}$, consider $\varepsilon'\downarrow0$, and notice that $\Lambda_{\sigma^\eta+\varepsilon'}\rightarrow\Lambda_{\sigma^\eta}$ and that the integrals in the last expression of the above vanish as $\varepsilon'\downarrow0$, to obtain $\Lambda_{\sigma^\eta}^{d-1}\leq0$.~This is the desired contradiction.
\qed

\subsection{Lower bound on the jump sizes}\label{subse:lowerBound}

In this subsection, we show \eqref{jump_cond} for any solution $(\Lambda,u)$ to \eqref{what_is_u}, \eqref{growth} with $d\geq3$ satisfying \eqref{eq.Sec4.DeltaLambda.upper}.

\begin{proposition}\label{prop_LBD}
Let $(\Lambda,u)$ be a solution to \eqref{what_is_u}, \eqref{growth} with $d\geq3$ satisfying \eqref{eq.Sec4.DeltaLambda.upper}, and let $t\in[0,T\wedge\zeta)$ be a jump time of $\Lambda$.~Then,
\begin{equation}
\Lambda_{t-}-\Lambda_t \ge \inf\bigg\{y\in(0,\Lambda_{t-}):\, \int_{\Lambda_{t-}-y}^{\Lambda_{t-}} u(t-,x)\,\nu(\mathrm{d}x) >\int_{\Lambda_{t-}-y}^{\Lambda_{t-}} \big(H(x)-1\big)\,\nu(\mathrm{d}x) \bigg\}.
\end{equation}
\end{proposition}

\smallskip

\noindent\textbf{Proof.}
Recall that Lemmas \ref{le:Sec4.u.bound.H} and \ref{le:noJumps.smallU} imply that $\Lambda$ cannot have jumps before $\sigma$. Therefore, $t\geq\sigma$. Recall also that $\Lambda$ is non-increasing on $[\sigma-,T\wedge\zeta)$. Hence, $\Lambda_{t-}\!>\!\Lambda_t$ and $u(s,\cdot)\!\leq\! H(\Lambda_{s})$ for all sufficiently small $s\geq t$.

\smallskip

Then, choosing $r<t<s$ sufficiently close to $t$ and using \eqref{eq.Sec4.FK.t.s}, \eqref{growth} we obtain
\begin{align}
&\;\frac{1}{d}\big((\Lambda_r)^d-(\Lambda_s)^d\big) = \int_{\RR_+} u(s,x)\,\nu(\mathrm{d}x) 
- \int_{\RR_+} u(r,x)\,\nu(\mathrm{d}x) \nonumber \\
& = \int_{\RR_+} \E^x\big[\mathbf{1}_{\{\tau_{\Lambda_{s-\cdot}}\leq s-r\}}\,H(R_{\tau_{\Lambda_{s-\cdot}}}) + \mathbf{1}_{\{\tau_{\Lambda_{s-\cdot}}>s-r\}}\,u(r,R_{s-r}) - u(r,x) \big]\,\nu(\mathrm{d}x) \nonumber \\
& = \int_{\RR_+} \E^x\big[\mathbf{1}_{\{\tau_{\Lambda_{s-\cdot}}\leq s-r\}}\,(H(R_{\tau_{\Lambda_{s-\cdot}}})-u(r,R_{s-r}))
+ u(r,R_{s-r}) - u(r,x) \big]\,\nu(\mathrm{d}x) \nonumber \\
& = \int_{\RR_+} \E^x\big[\mathbf{1}_{\{\tau_{\Lambda_{s-\cdot}}\leq s-r\}}\,(H(R_{\tau_{\Lambda_{s-\cdot}}})-u(r,R_{s-r}))\big]\,\nu(\mathrm{d}x)\label{eq.JumpLowerBnd.eq2} \\
& = \frac{\Gamma(d/2+1)}{d\pi^{d/2}} \int_{\RR^d} \E\big[\mathbf{1}_{\{\tau_s(x)\leq s-r\}}\,\big(H(|x+B_{\tau_s(x)}|)-u(r,|x+B_{s-r}|)\big)\big]\,\mathrm{d}x \nonumber \\
& = \frac{\Gamma(d/2+1)}{d\pi^{d/2}} \int_{\RR^d} \E\big[\mathbf{1}_{\{\tau_s(x-B_{s-r})\leq s-r\}}\,\big(H(|x+B_{\tau_s(x-B_{s-r})}-B_{s-r}|)-u(r,|x|)\big)\big]\,\mathrm{d}x, \nonumber
\end{align}
where $\Gamma$ is the Gamma function, $B$ is a standard Brownian motion in $\RR^d$, and $\tau_s(x)$ is the first hitting time of the boundary $\Lambda_{s-\cdot}$ by $|x+B_\cdot|$.~In the above, we have used Fubini's Theorem, which is justified by the arguments below.

\medskip

Notice that, almost surely for Lebesgue almost every $x\in\RR^d$, $u(r,|x|)$ converges to $u(t-,|x|)$ (recall Lemma \ref{le:u.LRlims}) and $\mathbf{1}_{\{\tau_s(x-B_{s-r})\leq s-r\}}$ converges to $\mathbf{1}_{(\Lambda_t,\Lambda_{t-})}(|x|)$, as $r\uparrow t$ and $s\downarrow t$.
In addition, as $\Lambda$ is bounded away from zero in $[-1,t+\varepsilon]$ for some $\varepsilon>0$, and in view of Lemma \ref{le:Sec4.u.bound.H}, we conclude that $H(|x+B_{\tau_s(x-B_{s-r})}-B_{s-r}|)-u(r,|x|)$ is bounded in absolute value uniformly over all $x>0$ and all small enough $s-r>0$. Moreover, $\PP(\tau_s(x-B_{s-r})\leq s-r)$ is bounded from above by an integrable function of $x$, uniformly over all small enough $s-r>0$.
Then, sending $s\downarrow t$, $r\uparrow t$, and using the Dominated Convergence Theorem we deduce from \eqref{eq.JumpLowerBnd.eq2}:
\begin{align}
\frac{1}{d}\big((\Lambda_{t-})^d\!-\!(\Lambda_t)^d\big) 
\!=\! \frac{\Gamma(d/2\!+\!1)}{d\pi^{d/2}} \int_{\{\Lambda_t<|x|<\Lambda_{t-}\}} \! H(|x|) \!-\! u(t-,|x|)\,\mathrm{d}x
\!=\! \int_{\Lambda_t}^{\Lambda_{t-}} \! H(x) \!-\! u(t-,x)\,\nu(\mathrm{d}x).\label{eq.JumpLowerBnd.eq10n}
\end{align}

\smallskip

Now, let us assume that there exists an $\eta\in(0,\Lambda_{t-}]$ such that $0<\Lambda_{t-}-\Lambda_{t} < \eta$ and
\begin{align}\label{eq.JumpLowerBnd.eq1}
\int_{\Lambda_{t-}-y}^{\Lambda_{t-}} u(t-,x)\,\nu(\mathrm{d}x) \leq \int_{\Lambda_{t-}-y}^{\Lambda_{t-}} \big(H(x)-1\big)\,\nu(\mathrm{d}x),\quad y\in[0,\eta).
\end{align}
Our goal is to show that this assumption leads to a contradiction, yielding the proposition. To this end, we observe that \eqref{eq.JumpLowerBnd.eq10n} gives
\begin{align}\label{dens_eq}
\int_{\Lambda_{t}}^{\Lambda_{t-}} u(t-,x)\,\nu(\mathrm{d}x) = \int_{\Lambda_{t}}^{\Lambda_{t-}} \big(H(x)-1\big)\,\nu(\mathrm{d}x).
\end{align}
Let us now fix a $\chi\in(0,\eta-(\Lambda_{t-}-\Lambda_{t}))$ and recall that $u(t-,x)=u(t,x)$ for $x<\Lambda_t$ (see Lemma \ref{le:u.LRlims}), which along with \eqref{dens_eq} and \eqref{eq.JumpLowerBnd.eq1} implies
\begin{align}\label{eq.JumpLowerBnd.eq1n}
\int_{\Lambda_{t}-y}^{\Lambda_{t}} u(t,x)\,\nu(\mathrm{d}x) \leq \int_{\Lambda_{t}-y}^{\Lambda_{t}} \big(H(x)-1\big)\,\nu(\mathrm{d}x),
\quad y\in[0,\chi].
\end{align}

\smallskip

Next, we repeat the first four steps in \eqref{eq.JumpLowerBnd.eq2} with $r=t<s$ and use the monotonicity of $\Lambda$ on $[t,s]$ to obtain
\begin{align*}
& \frac{1}{d}\big((\Lambda_t)^d-(\Lambda_s)^d\big) 
\geq  \int_{\Lambda_s}^{\infty} \E^x\big[\mathbf{1}_{\{\tau_{\Lambda_{s-\cdot}}\leq s-t\}}\,(H(R_{\tau_{\Lambda_{s-\cdot}}})-u(t,R_{s-t}))\big]\,\nu(\mathrm{d}x).
\end{align*}
Due to the right-continuity of $\Lambda$, for $x>\Lambda_s$ the event $R_{s-t}\in [\Lambda_{t}-\chi,\Lambda_t]$ implies $\tau_{\Lambda_{s-\cdot}}\leq s-t$. Thus, for small enough $s-t>0$ and $\chi<\Lambda_{t-}-\Lambda_t$ we have
\begin{align*}
\frac{1}{d}\big((\Lambda_t)^d-(\Lambda_s)^d\big) 
& \ge \int_{\Lambda_s}^{\infty} \E^x\big[\mathbf{1}_{\{R_{s-t}\in [\Lambda_{t}-\chi,\Lambda_t]\}}\,(H(R_{s-t})-u(t,R_{s-t}))\big]\,\nu(\mathrm{d}x) \\
& \;\;\;\, + \int_{\Lambda_s}^{\infty} \E^x\big[\mathbf{1}_{\{\tau_{\Lambda_{s-\cdot}}\leq s-t\}}\,(H(R_{\tau_{\Lambda_{s-\cdot}}})-H(R_{s-t}))\big]\,\nu(\mathrm{d}x).
\end{align*}
Notice also that $H(R)$ is a supermartingale for $d\geq3$. Hence,
\begin{align*}
&\;\E^x\big[\mathbf{1}_{\{\tau_{\Lambda_{s-\cdot}}\leq s-t\}}\,
(H(R_{\tau_{\Lambda_{s-\cdot}}})-H(R_{s-t}))\big] \\
& = \E^x\big[\mathbf{1}_{\{\tau_{\Lambda_{s-\cdot}}\leq s-t\}}\,\EE^x[H(R_{\tau_{\Lambda_{s-\cdot}}\wedge (s-t)})-H(R_{s-t})
\,|\,\mathcal{F}^R_{\tau_{\Lambda_{s-\cdot}}}]\big]
\geq0,
\end{align*}
and so,
\begin{align}
\frac{1}{d}\big((\Lambda_t)^d-(\Lambda_s)^d\big) 
\geq  \int_{\Lambda_s}^{\infty} \E^x\big[\mathbf{1}_{\{R_{s-t}\in [\Lambda_{t}-\chi,\Lambda_t]\}}\,(H(R_{s-t})-u(t,R_{s-t}))\big]\,\nu(\mathrm{d}x).\label{eq.JumpLowerBnd.t.to.s.dropSupermtg}
\end{align}

\smallskip

Further, with $f(z):=\mathbf{1}_{\{|z|\in [\Lambda_t-\chi,\Lambda_t]\}}(H(|z|)-u(t,|z|))$, $z\in\RR^d$, we have for $\Lambda_t-\Lambda_s<\chi/8$:
\begin{align*}
 \int_{\Lambda_s}^{\infty} \E^x\big[\mathbf{1}_{\{R_{s-t}\in [\Lambda_{t}-\chi,\Lambda_t]\}}\,(H(R_{s-t})-u(t,R_{s-t}))\big]\,\nu(\mathrm{d}x)
= \frac{\Gamma(d/2\!+\!1)}{d\pi^{d/2}}\int_{\{|z|\geq\Lambda_s\}} \E[f(z\!+\!B_{s-t})]\,\mathrm{d}z \\
= \frac{\Gamma(d/2+1)}{d\pi^{d/2}}\,\E\bigg[\int_{\RR^d} \bone_{\{|z-B_{s-t}| \geq \Lambda_s\}} \mathbf{1}_{\{|z|\in [\Lambda_{t}-\chi,\Lambda_t]\}} (H(|z|)-u(t,|z|))\,\mathrm{d}z\bigg] \\
\geq \frac{\Gamma(d/2+1)}{d\pi^{d/2}}\, \E\bigg[\bone_{\{|B_{s-t}|\leq \chi/2\}} \int_{\RR^d} \bone_{\{|z-B_{s-t}| \geq \Lambda_s,\,|z|\leq \Lambda_t\}} (H(|z|)-u(t,|z|))\,\mathrm{d}z\bigg] \\
= \frac{\Gamma(d/2+1)}{d\pi^{d/2}}\,\E\bigg[\bone_{\{|B_{s-t}|\leq \chi/2\}} \int_{S^{d-1}} \int_{\phi(v)\wedge\Lambda_t}^{\Lambda_t}  \left(H(r)-u(t,r)\right)\,r^{d-1}\,\mathrm{d}r\,\sigma(\mathrm{d}v)\bigg],
\end{align*}
where $B$ is a standard $d$-dimensional Brownian motion, $S^{d-1}$ is the $(d-1)$-dimensional unit sphere centered at the origin, $\sigma(\mathrm{d}v)$ is the volume element of this sphere, and $\phi(v)$ is the distance between the origin and the sphere $\{z\!:|z-B_{s-t}| = \Lambda_s\}$ along the direction $v\in S^{d-1}$. On the event $\{|B_{s-t}|\leq \chi/2\}$, we have $\Lambda_t-\phi(v)\leq\chi$, and hence \eqref{eq.JumpLowerBnd.eq1n} yields
\begin{align*}
& \int_{S^{d-1}} \int_{\phi(v)\wedge\Lambda_t}^{\Lambda_t}  (H(r)-u(t,r))\,r^{d-1}\,\mathrm{d}r\,\sigma(\mathrm{d}v)
\geq \int_{S^{d-1}} \int_{\phi(v)\wedge\Lambda_t}^{\Lambda_t} r^{d-1}\mathrm{d}r\,\sigma(\mathrm{d}v).
\end{align*}
The right-hand side is the Lebesgue measure of the difference $\{z\!:|z|\leq\Lambda_t\}\setminus\{z\!:|z\!-\!B_{s-t}|\!\geq\!\Lambda_s\}$ between two balls, which is equal to $\frac{\pi^{d/2}}{\Gamma(d/2+1)}((\Lambda_t)^d - (\Lambda_s)^d)$, provided $|B_{s-t}|\leq \Lambda_t-\Lambda_s$. Whenever the latter inequality is violated, the aforementioned measure is strictly larger than $\frac{\pi^{d/2}}{\Gamma(d/2+1)}((\Lambda_t)^d - (\Lambda_s)^d)$, and it depends continuously on $|B_{s-t}|$.~Using this observation and recalling that $\Lambda_t-\Lambda_s<\chi/8$ we obtain
\begin{align*}
&\,\frac{\Gamma(d/2+1)}{d\pi^{d/2}}\,\E\bigg[\bone_{\{|B_{s-t}|\leq\chi/2\}} \int_{S^{d-1}} \int_{\phi(v)\wedge\Lambda_t}^{\Lambda_t}  \left(H(r)-u(t,r)\right)\,r^{d-1}\,\mathrm{d}r\,\sigma(\mathrm{d}v)\bigg] \\
&\geq \frac{1}{d}\big((\Lambda_t)^d-(\Lambda_s)^d\big)\,\PP(|B_{s-t}|\leq \chi/2)
+ C_1\PP(|B_{s-t}|\in [\chi/4, \chi/2]) \\
&\geq \frac{1}{d}\big((\Lambda_t)^d-(\Lambda_s)^d\big)
+ C_1\PP(|B_{s-t}|>\chi/4) - C_2 \PP(|B_{s-t}|>\chi/2),
\end{align*}
with some $C_1,C_2\in(0,\infty)$.~For sufficiently small $s-t>0$, the above is strictly larger than $\frac{1}{d} ((\Lambda_t)^d-(\Lambda_s)^d)$, which yields the desired contradiction to \eqref{eq.JumpLowerBnd.t.to.s.dropSupermtg}.
\qed

\subsection{Upper bound on the jump sizes}\label{subse:upperBound}

In this subsection, we show that the absolute jump sizes of the limits points of $\{\Lambda^\Delta\}_{\Delta\downarrow0}$ satisfy the desired upper bound.~The following proposition is the main result of this subsection and completes the proof of Theorem \ref{thm_main}.

\begin{proposition}\label{prop_UBD}
Consider a limit point $\Lambda$ of $\{\Lambda^\Delta\}_{\Delta\downarrow0}$, a jump time $t\in[0,T\wedge\zeta]$ of $\Lambda$, and $u$ defined by \eqref{what_is_u}. Then,
\begin{align}
&\Lambda_{t-}-\Lambda_t \le  \inf\bigg\{y\in(0,\Lambda_{t-}]:\, \int_{\Lambda_{t-}-y}^{\Lambda_{t-}} u(t-,x)\,\nu(\mathrm{d}x) >\int_{\Lambda_{t-}-y}^{\Lambda_{t-}} \big(H(x)-1\big)\,\nu(\mathrm{d}x) \bigg\},\label{eq.JmpUpperBnd.mainProp}\\
& \Lambda_{t}-\Lambda_{t-} \leq \inf\bigg\{y>0:\, \int_{\Lambda_{t-}}^{\Lambda_{t-}+y} u(t-,x)\,\nu(\mathrm{d}x) <\int_{\Lambda_{t-}}^{\Lambda_{t-}+y} \big(H(x)+1\big)\,\nu(\mathrm{d}x)\bigg\}. \label{eq.JmpUpperBnd.mainProp.n}
\end{align}
\end{proposition}


\smallskip

The proof of Proposition \ref{prop_UBD} relies on the following lemmas.

\begin{lemma}\label{down_jumps_small}
There are $h_1,h_2\!:(0,\infty)\to\RR_+$ bounded on compact subsets of $(0,\infty)$ such that
\begin{equation}\label{Lemma4.5ineq}
\begin{split}
\int_{\RR_+} u^\Delta((m+1)\Delta-,x)\,\nu(\mathrm{d}x)
- \int_{\RR_+} u^\Delta(m\Delta,x)\,\nu(\mathrm{d}x)
\le h_1(\Lambda^\Delta_{m\Delta})\,\sqrt{\Delta},\\
m=0,\,1,\,\ldots,\;\;\Delta\in\big(0,1/h_2(\Lambda^\Delta_{m\Delta})\big)\;\;\text{such that} \;\;\Lambda^\Delta_{m\Delta}>0.
\end{split}
\end{equation} 
\end{lemma}

\smallskip

\noindent\textbf{Proof.}
By using Definition \ref{def:Euler}, we represent the left-hand side in \eqref{Lemma4.5ineq} as
\begin{equation*}
	\begin{split}
		&\; \int_{\RR_+} \E^x\big[\mathbf{1}_{\{\tau_{\Lambda^\Delta_{m\Delta}}<\Delta\}}\big(H(\Lambda^\Delta_{m\Delta})-u^\Delta(m\Delta,R_\Delta)\big)\big]\,\nu(\mathrm{d}x) \\
		& + \int_{\RR_+} \E^x[u^\Delta(m\Delta,R_\Delta)]\,\nu(\mathrm{d}x)
		-\int_{\RR_+} u^\Delta(m\Delta,x)\,\nu(\mathrm{d}x).
	\end{split}
\end{equation*}
Since the $d$-dimensional Bessel process is reversible with respect to the measure $\nu$, the second line in the latter display is equal to $0$. Therefore, we arrive at the upper bound
\begin{equation}\label{hitting_time_bnd}
	H(\Lambda^\Delta_{m\Delta})\,
	\int_{\RR_+} \PP^x(\tau_{\Lambda^\Delta_{m\Delta}}<\Delta)\,\nu(\mathrm{d}x).
\end{equation}

\smallskip

Finally, we split the integral in \eqref{hitting_time_bnd} according to $x\ge\Lambda^\Delta_{m\Delta}$, $x\in[\sqrt{(\Lambda^\Delta_{m\Delta})^2-d\Delta},\Lambda^\Delta_{m\Delta})$, and $x\in[0,\sqrt{(\Lambda^\Delta_{m\Delta})^2-d\Delta})$. For $x\ge\Lambda^\Delta_{m\Delta}$, we bound $\PP^x(\tau_{\Lambda^\Delta_{m\Delta}}<\Delta)$ by replacing the Bessel process $R^x$ with a standard Brownian motion started from $x$. For $x\in[\sqrt{(\Lambda^\Delta_{m\Delta})^2-d\Delta},\Lambda^\Delta_{m\Delta})$, we estimate $\PP^x(\tau_{\Lambda^\Delta_{m\Delta}}<\Delta)$ by $1$. For $x\in[0,\sqrt{(\Lambda^\Delta_{m\Delta})^2-d\Delta})$, we apply the Dambis-Dubins-Schwarz Theorem (see, e.g., \cite[Chapter 3, Problem 4.7]{KaSh}) upon noting that the diffusion coefficient of $(R^x)^2$ is smaller or equal to $4(\Lambda^\Delta_{m\Delta})^2$ until $\tau_{\Lambda^\Delta_{m\Delta}}$, and that its drift coefficient equals to $d$. Thus, we get the upper bound
\begin{equation*}
	\begin{split}
		&\,\int_{\Lambda^\Delta_{m\Delta}}^\infty\!\! 2\overline{\Phi}\bigg(\frac{x\!-\!\Lambda^\Delta_{m\Delta}}{\sqrt{\Delta}}\bigg)\,\nu(\mathrm{d}x)
		\!+\!\int_{\sqrt{(\Lambda^\Delta_{m\Delta})^2-d\Delta}}^{\Lambda^\Delta_{m\Delta}} \nu(\mathrm{d}x)
		\!+\!\int_0^{\sqrt{(\Lambda^\Delta_{m\Delta})^2-d\Delta}} \!2\overline{\Phi}\bigg(\frac{(\Lambda^\Delta_{m\Delta})^2\!-\!d\Delta\!-\!x^2}{\sqrt{\Delta}/(2\Lambda^\Delta_{m\Delta})}\bigg)\,\nu(\mathrm{d}x) \\
		& = 2\sqrt{\Delta}\int_0^\infty \overline{\Phi}(y)(\Lambda^\Delta_{m\Delta}+y\sqrt{\Delta})^{d-1}\,\mathrm{d}y
		+\frac{1}{d}\big((\Lambda^\Delta_{m\Delta})^d-((\Lambda^\Delta_{m\Delta})^2-d\Delta)^{d/2}\big) \\
		&\quad
		+\frac{\sqrt{\Delta}}{2\Lambda^\Delta_{m\Delta}}\int_0^{\frac{(\Lambda^\Delta_{m\Delta})^2-d\Delta}{\sqrt{\Delta}/(2\Lambda^\Delta_{m\Delta})}} \overline{\Phi}(y)\,\bigg((\Lambda^\Delta_{m\Delta})^2-d\Delta-\frac{y\sqrt{\Delta}}{2\Lambda^\Delta_{m\Delta}}\bigg)^{d/2-1}\,\mathrm{d}y.
	\end{split}
\end{equation*}
The latter expression readily admits an estimate of the form asserted in \eqref{Lemma4.5ineq}. \qed 

\begin{lemma}\label{down_jumps_small.2}
For any $\varepsilon>0$, there exist $C=C(\varepsilon)$, $\Delta_0=\Delta_0(\varepsilon)$, both in $(0,\infty)$, such that
\begin{align*}
&\Lambda^\Delta_{m\Delta} - \Lambda^\Delta_{(m+1)\Delta}
\le C\sqrt{\Delta},
\end{align*} 
provided $(\Lambda^\Delta_{m\Delta}/\gamma)\in[1+\varepsilon,\,1/\varepsilon]$ and $\Delta\in(0,\Delta_0)$.
\end{lemma}

\smallskip

\noindent\textbf{Proof.} By Definition \ref{def:Euler}, $\Lambda^\Delta_{(m+1)\Delta}$ is bounded from below by any $y\in [0,\Lambda^\Delta_{m\Delta})$ satisfying
\begin{align*}
& \,\int_{\RR_+} u^\Delta((m+1)\Delta-,x)-u^\Delta(m\Delta,x)\,\nu(\mathrm{d}x)
+\int_{y}^{\Lambda^\Delta_{m\Delta}} H(x)-u^\Delta((m+1)\Delta-,x)\,\nu(\mathrm{d}x) \\
& < \frac{1}{d}\big((\Lambda^\Delta_{m\Delta})^d-y^d\big).
\end{align*}
Lemma \ref{down_jumps_small} and the fact that $H(y)\!\leq\! 1/(1\!+\!\varepsilon/2)$ for $y\in [\Lambda^\Delta_{m\Delta}\!-\!\varepsilon\gamma/2,\Lambda^\Delta_{m\Delta}]$ yields for such $y$:
\begin{align*}
& \int_{\RR_+} u^\Delta((m+1)\Delta-,x)-u^\Delta(m\Delta,x)\,\nu(\mathrm{d}x)
+\int_{y}^{\Lambda^\Delta_{m\Delta}} H(x)-u^\Delta((m+1)\Delta-,x)\,\nu(\mathrm{d}x) \\
& \leq C_1(\varepsilon) \sqrt{\Delta} + \int_{y}^{\Lambda^\Delta_{m\Delta}} 1/(1+\varepsilon/2)\,\nu(\mathrm{d}x)
= C_1(\varepsilon) \sqrt{\Delta} + \frac{1}{d(1+\varepsilon/2)}\big((\Lambda^\Delta_{m\Delta})^d\!-\!y^d\big).
\end{align*}
Thus, 
\begin{align*}
& (\Lambda^\Delta_{(m+1)\Delta})^d \geq (\Lambda^\Delta_{m\Delta})^d - \frac{d(1+\varepsilon/2) C_1(\varepsilon)}{\varepsilon/4} \sqrt{\Delta},
\end{align*}
which implies, for all $\Delta\in(0,\Delta_0(\varepsilon))$,
\begin{align*}
& \qquad\qquad\;\, d(\gamma+\gamma\varepsilon/2)^{d-1}(\Lambda^\Delta_{m\Delta} - \Lambda^\Delta_{(m+1)\Delta}) \leq (\Lambda^\Delta_{m\Delta})^d - (\Lambda^\Delta_{(m+1)\Delta})^d \leq C_2(\varepsilon) \sqrt{\Delta}.
\qquad\qquad\;\,\qed
\end{align*}

\smallskip

\begin{lemma}\label{down_jumps_small.3}
For any $t\in[\sigma,T\wedge\zeta]$ (with $\sigma$ defined in \eqref{eq.Sec4.sigma.def}) and $n\in\NN$, let 
\begin{equation*}
z_n:=\inf\bigg\{y\in(0,\Lambda_{t-}]:\,\int_{\Lambda_{t-}-y}^{\Lambda_{t-}} u(t-,x)\,\nu(\mathrm{d}x) >\int_{\Lambda_{t-}-y}^{\Lambda_{t-}} \big(H(x)-1\big)\,\nu(\mathrm{d}x) + \frac{1}{n}\bigg\},
\end{equation*}
and, for $r\in[-1,t)$, 
\begin{align*}
& \widetilde\sigma=\widetilde\sigma(\Delta,r,n):=\inf\{s\in[r,T\wedge\zeta]:\, \Lambda^\Delta_s\leq \Lambda_{t-}-z_n\}.
\end{align*}
Then, for any $n\in\nn$ large enough, we have
\begin{align*}
\underline{\lim}_{r\uparrow t}\,
\underline{\lim}_{\Delta\downarrow0}\Lambda^\Delta_{\widetilde\sigma}\geq \Lambda_{t-}-z_n
\end{align*}
if $\Lambda_{t-}-\Lambda_t>z_n$, where the first limit inferior is taken over the continuity points $r$ of $\Lambda$.
\end{lemma}

\smallskip

\noindent\textbf{Proof.}
Recall the condition that determines the size of a downward jump:
\begin{equation}\label{eq.JmpUpperBd.le3.EulJump.def} 
\begin{split}
\Lambda^\Delta_{\widetilde\sigma}=
0\vee\sup\bigg\{y\in[0,\Lambda^\Delta_{\widetilde\sigma-\Delta}):\;\int_{\RR_+} &\; u^\Delta(\widetilde\sigma-,x)\,\mathbf{1}_{\RR_+\backslash [y,\Lambda^\Delta_{\widetilde\sigma-\Delta}]}(x)
+H(x)\,\mathbf{1}_{[y,\Lambda^\Delta_{\widetilde\sigma-\Delta}]}(x) \\
& -u^\Delta(\widetilde\sigma-\Delta,x)\,\nu(\mathrm{d}x)
<\frac{1}{d}\big((\Lambda^\Delta_{\widetilde\sigma-\Delta})^d-y^d\big)\bigg\}.
\end{split}
\end{equation}
Due to Lemma \ref{down_jumps_small}, for any $y\in[0,\Lambda^\Delta_{\widetilde\sigma-\Delta})$,
\begin{equation}\label{eq.JmpUpperBd.le3.eq0}
\begin{split}
&\;\int_{\RR_+} u^\Delta(\widetilde\sigma-,x)\,\mathbf{1}_{\RR_+\backslash [y,\Lambda^\Delta_{\widetilde\sigma-\Delta}]}(x)
+H(x)\,\mathbf{1}_{[y,\Lambda^\Delta_{\widetilde\sigma-\Delta}]}(x) - u^\Delta(\widetilde\sigma-\Delta,x)\,\nu(\mathrm{d}x) \\
&=\int_{\RR_+} u^\Delta(\widetilde\sigma-,x) - u^\Delta(\widetilde\sigma-\Delta,x)\,\nu(\mathrm{d}x)
+\int_{\RR_+}(H(x)-u^\Delta(\widetilde\sigma-,x))\,\mathbf{1}_{[y,\Lambda^\Delta_{\widetilde\sigma-\Delta}]}(x)\,\nu(\mathrm{d}x) \\
& \leq C_1 \sqrt{\Delta}
+\int_{y}^{\Lambda^\Delta_{\widetilde\sigma-\Delta}} H(x)-u^\Delta(\widetilde\sigma-,x)\,\nu(\mathrm{d}x).
\end{split}
\end{equation}

\smallskip

We argue by contradiction and assume that there exist an $\varepsilon_0\in(0,\Lambda_{t-}-z_n]$ and a sequence $r_m\uparrow t$ of continuity points of $\Lambda$ such that for every $r=r_m$,
\begin{align}
\underline{\lim}_{\Delta\downarrow0}\Lambda^\Delta_{\widetilde\sigma}\leq \Lambda_{t-}-z_n - \varepsilon_0. \label{eq.Sec4.Le4.7.contra}
\end{align}
Then, we choose $m_0\in\NN$ large enough to ensure sure that $\inf_{s\in[r_{m_0},t)}\Lambda_s > \Lambda_{t-}-z_n$ (note that $z_n>0$).~As a consequence, there exists a $\Delta_0(m_0)$ so that for all $m\geq m_0$ and all $\Delta\in(0,\Delta_0(m_0))$, it holds $\widetilde\sigma\in(r_m,T\wedge\zeta]$, and hence $\Lambda^\Delta_{\widetilde\sigma-\Delta}>\Lambda_{t-}-z_n$.~Next, we use the definition of $z_n$ to deduce the existence of $0<\delta_1<\delta_2<\varepsilon_0$ and $\delta_3>0$ with the property
\begin{equation*}
\inf_{y'\in[\Lambda_{t-}-z_n,\,\Lambda_{t-}],\,y\in [\Lambda_{t-}-z_n-\delta_2,\,\Lambda_{t-}-z_n-\delta_1]}\int_{y}^{y'} u(t-,x)\,-\,\big(H(x)-1\big)\,\nu(\mathrm{d}x)\geq 2\delta_3.
\end{equation*}
Let us fix an arbitrary $y\in [\Lambda_{t-}-z_n-\delta_2,\Lambda_{t-}-z_n-\delta_1]$.~Increasing $m_0$ and making it depend on $y$, if needed, we combine the convergence of $u(r_m,\cdot)$ to $u(t-,\cdot)$ (see Lemma \ref{le:u.LRlims}) and the Dominated Convergence Theorem (recall Lemma \ref{le:Sec4.u.bound.H}) to conclude that
\begin{align*}
&\inf_{y'\in[\Lambda_{t-}-z_n,\Lambda_{t-}]}\int_{y}^{y'} u(r_m,x)\,-\,\big(H(x)-1\big)\,\nu(\mathrm{d}x) \geq \delta_3,\quad m\geq m_0.
\end{align*}

Now, from the M1-convergence $\Lambda^\Delta\!\to\!\Lambda$ we infer that for all $m\!\geq\! m_0$, we have~$\lim_{\Delta\downarrow0}\widetilde\sigma\!=\!t$. This observation, as well as the M1-convergence of $\Lambda^\Delta\to\Lambda$, the property $\Lambda^\Delta_{\widetilde\sigma-\Delta}>\Lambda_{t-}-z_n$, and the boundedness of $u$ on $[0,t)\times \RR_+$ (see Lemma \ref{le:Sec4.u.bound.H}) and of $H$ on $[\Lambda_{t-}-z_n-\delta_2,\infty)$ yield via the latter display: 
\begin{equation}\label{eq.JmpUpperBd.le3.eq3}
\underline{\lim}_{\Delta\downarrow0} \int_{y}^{\Lambda^\Delta_{\widetilde\sigma-\Delta}} u(r_m,x)\,-\,\big(H(x)-1\big)\,\nu(\mathrm{d}x)\geq\delta_3,\quad m\ge m_0. 
\end{equation}
We consider two possible cases.~The first case is $t=\sigma$.~In this case, Lemma \ref{le:Sec4.u.bound.H} and the choice of $z_n$ imply $\Lambda_{t-}-z_n>\gamma$ for large enough $n\geq 1$ (as $\Lambda_{t-}>\gamma$ and $\lim_{n\rightarrow\infty}z_n=0$). Then, for any sufficiently small $\varepsilon>0$ and any $m\geq m_0$, there exists a $\Delta_1(m)>0$ such that $\Lambda^\Delta_{\widetilde\sigma-}=\Lambda^\Delta_{\widetilde\sigma-\Delta}\geq \Lambda_{t-}-z_n \geq \gamma(1+\varepsilon)$ for all $\Delta\in(0,\Delta_1(m))$.~Finally, Lemma \ref{down_jumps_small.2} and the inequality $\Lambda^\Delta_{\widetilde\sigma-}\geq \Lambda_{t-}-z_n$ give the desired contradiction to \eqref{eq.Sec4.Le4.7.contra}.

\medskip

The second case is $t>\sigma$.~
Increasing $m_0$ we ensure that $r_{m_0}>\sigma$.~Decreasing $\Delta_0>0$, if needed, we deduce from Lemma \ref{le:osc} that, for any $m\geq m_0$ and all $\Delta\in(0,\Delta_0)$, the function~$\Lambda^\Delta$ is non-increasing on $[r_{m_0},\infty)$ and $u^\Delta(r_m,\cdot)\leq H(\Lambda^\Delta_{r_m})$.~Using this observation and the Feynman-Kac representation we obtain, for any $x>0$, $m\geq m_0$, $\Delta\in(0,\Delta_0)$, and $\xi\in(0,(t-r_{m})\wedge\Delta)$:
\begin{align*}
u^\Delta(\widetilde\sigma-\xi,x) &= \E^x\big[\mathbf{1}_{\{\tau_{\Lambda^\Delta_{\widetilde\sigma-\xi-\cdot}}\leq\widetilde\sigma-\xi-r_m\}}
\,H(R_{\tau_{\Lambda^\Delta_{\widetilde\sigma-\xi-\cdot}}})\big]
+\E^x\big[\mathbf{1}_{\{\tau_{\Lambda^\Delta_{\widetilde\sigma-\xi-\cdot}}> \widetilde\sigma-\xi-r_m\}}\,u^\Delta(r_m,R_{\widetilde\sigma-\xi-r_m})\big]\\
&\geq \E^x[u^\Delta(r_m,R_{\widetilde\sigma-\xi-r_m})].
\end{align*}
Sending $\xi\downarrow0$ in the latter display and using the continuity of the transition kernel of $R$ in time with respect to the $L^1$-norm in space together with the boundedness of $u^\Delta(r_m,\cdot)$ for sufficiently small $\Delta$ (see Remark \ref{rem:uDelta.bound}) we deduce that for any $x>0$, $m\geq m_0$, $\Delta\in(0,\Delta_1(m))$,
\begin{align*}
&u^\Delta(\widetilde\sigma-,x) \geq \E^x[u^\Delta(r_m,R_{\widetilde\sigma-r_m})].
\end{align*}

\smallskip

Next, we take $\Delta\downarrow0$ and rely on $\widetilde\sigma\rightarrow t$ and $u^\Delta(r_m,\cdot)\rightarrow u(r_m,\cdot)$ (pointwise in $\RR_+\setminus\{\Lambda_{r_m}\}$, cf.~Lemma \ref{lem:u_conv}), as well as the Dominated Convergence Theorem and the continuity of the transition kernel of $R$ in time to conclude that for any $m\geq m_0$ and $\Delta\in(0,\Delta_1(m))$,
\begin{align*}
\int_y^{\Lambda^\Delta_{\widetilde\sigma-\Delta}} \big|\E^x[u^\Delta(r_m,R_{\widetilde\sigma-r_m})]
-\E^x[u(r_m,R_{t-r_m})]\big|\,\nu(\mathrm{d}x)\le\frac{\delta_3}{4},
\end{align*}
where we decrease $\Delta_1(m)>0$, if needed.~Collecting the above, we obtain
\begin{align}\label{eq.JmpUpperBd.le3.eq1}
& \int_y^{\Lambda^\Delta_{\widetilde\sigma-\Delta}} u^\Delta(\widetilde\sigma-,x)\,\nu(\mathrm{d}x) 
\geq
\int_y^{\Lambda^\Delta_{\widetilde\sigma-\Delta}} \E^x[u(r_m,R_{t-r_m})]\,\nu(\mathrm{d}x)- \frac{\delta_3}{4},
\end{align}
for any $m\geq m_0$ and $\Delta\in(0,\Delta_1(m))$.

\medskip

Now, we take the functions $\widehat u(r_m,z)\!:=\!u(r_m,|z|)$ on $\RR^d$, and increase $m_0$ and decrease~$\Delta_1(m)$, if needed, to ensure that
\begin{equation}\label{eq.JmpUpperBd.le3.eq2}
\begin{split}
\int_y^{\Lambda^\Delta_{\widetilde\sigma-\Delta}} \E^x[u(r_m,R_{t-r_m})]\,\nu(\mathrm{d}x)
= \frac{\Gamma(d/2+1)}{d\pi^{d/2}} \int_{\{|z|\in[y,\Lambda^\Delta_{\widetilde\sigma-\Delta}]\}} \EE[\widehat u(r_m,z+B_{t-r_m})]\,\mathrm{d}z \\
=\frac{\Gamma(d/2+1)}{d\pi^{d/2}}\, \EE\bigg[\int_{\{|z-B_{t-r_m}|\in[y,\Lambda^\Delta_{\widetilde\sigma-\Delta}]\}} \widehat u(r_m,z)\,\mathrm{d}z\bigg] \\
\geq \frac{\Gamma(d/2+1)}{d\pi^{d/2}} \int_{\{|z|\in[y,\Lambda^\Delta_{\widetilde\sigma-\Delta}]\}} \widehat u(r_m,z)\,\mathrm{d}z-\frac{\delta_3}{4} 
= \int_y^{\Lambda^\Delta_{\widetilde\sigma-\Delta}} u(r_m,x)\,\nu(\mathrm{d}x) - \frac{\delta_3}{4}
\end{split}
\end{equation}
holds for all $m\geq m_0$ and $\Delta\in(0,\Delta_1(m))$.~In the above, we made use of the boundedness of $u$ on $[0,t)\times \RR_+$ and of the fact that the expected Lebesgue measure of the symmetric difference between the sets $\{z\!:|z-B_{t-r_m}|\in[y,\kappa]\}$ and $\{z\!:|z|\in[y,\kappa]\}$ vanishes as $r_m\uparrow t$, locally uniformly in $\kappa\geq y$.

\medskip

Collecting \eqref{eq.JmpUpperBd.le3.eq0}, \eqref{eq.JmpUpperBd.le3.eq1}, \eqref{eq.JmpUpperBd.le3.eq2}, \eqref{eq.JmpUpperBd.le3.eq3} we get that for any $y\in[\Lambda_{t-}-z_n-\delta_2,\Lambda_{t-}-z_n-\delta_1]$, there exists an $m_0>0$ such that for any $m\geq m_0$, there exists a $\Delta_1(m)>0$ for which
\begin{align*}
&\int_{\RR_+} u^\Delta(\widetilde\sigma-,x)\,\mathbf{1}_{\RR_+\backslash [y,\Lambda^\Delta_{\widetilde\sigma-\Delta}]}(x)
+H(x)\,\mathbf{1}_{[y,\Lambda^\Delta_{\widetilde\sigma-\Delta}]}(x) - u^\Delta(\widetilde\sigma-\Delta,x)\,\nu(\mathrm{d}x)\\
&\leq C_1 \sqrt{\Delta} + \int_{y}^{\Lambda^\Delta_{\widetilde\sigma-\Delta}} H(x)-u^\Delta(\widetilde\sigma-,x)\,\nu(\mathrm{d}x)\\
&\leq C_1 \sqrt{\Delta} + \int_{y}^{\Lambda^\Delta_{\widetilde\sigma-\Delta}} H(x)-\EE^x [u(r_m,R_{t-r_m})]\,\nu(\mathrm{d}x) + \frac{\delta_3}{4} \\
&\leq C_1 \sqrt{\Delta} + \int_{y}^{\Lambda^\Delta_{\widetilde\sigma-\Delta}} H(x)-u(r_m,x)\,\nu(\mathrm{d}x)
+ \frac{\delta_3}{2} \\
&\leq C_1 \sqrt{\Delta}
+\frac{1}{d}\big((\Lambda^\Delta_{\widetilde\sigma-\Delta})^d\!-\!y^d\big) - \frac{\delta_3}{4},
\quad \Delta\in(0,\Delta_1(m))
\end{align*}
holds. In view of \eqref{eq.JmpUpperBd.le3.EulJump.def}, this yields $\Lambda^\Delta_{\widetilde\sigma}\geq y\geq\Lambda_{t-}-z_n-\delta_2> \Lambda_{t-}-z_n-\varepsilon_0$, resulting in the desired contradiction to the choice of $\varepsilon_0$. \qed

\medskip

We conclude the subsection with the proof of Proposition \ref{prop_UBD}.

\medskip

\noindent\textbf{Proof of Proposition \ref{prop_UBD}.}
First, we recall that $\Lambda$ satisfies \eqref{growth}, with $u$ given by \eqref{what_is_u} (see Proposition \ref{prop:Stefan}).~Recall also (from Subsection \ref{subse:Sec4.prelim}) that $\Lambda$ does not jump in $[0,\sigma\wedge T)$, and that its jump at $\sigma$ (if any) must be downwards, whereby $\sigma$ is defined in \eqref{eq.Sec4.sigma.def}.~
On the other hand, Lemma \ref{le:osc} asserts that $\Lambda^\Delta$ is non-increasing in $[\sigma^\Delta,T]$ (see \eqref{eq.Sec2.sigmaDelta.def}).~The M1 convergence of $\Lambda^\Delta$ to $\Lambda$ and the above observations yield that $\Lambda$ is non-increasing on $(\sigma,T]$. Thus, to prove the proposition, it suffices to establish \eqref{eq.JmpUpperBnd.mainProp} for $t\in[\sigma,T\wedge\zeta]$.

\medskip

We let $z_n$, $n\in\NN$ be as in Lemma \ref{down_jumps_small.3}.~Notice that the sequence $(z_n)_{n\in\NN}$ is non-increasing and strictly positive.~If $z_n=\infty$ for all $n\in\NN$, then the right-hand side of \eqref{eq.JmpUpperBnd.mainProp} equals~$\infty$ and the statement of the proposition holds trivially.~Hence, without loss of generality we assume that $z_n\in(0,\Lambda_{t-}]$, $n\in\NN$,
and aim to show that $\Lambda_{t-}-\Lambda_t\le z_n$ for all $n\in\NN$.~For this purpose, we fix an $n\in\NN$ and suppose that $\Lambda_{t-}-\Lambda_t>z_n$.~To obtain a contradiction we pick any $\varepsilon\in(0,z_n)$ and any continuity time $0\le s<t$ of $\Lambda$ which is close enough to $t$ that $\Lambda_s>\Lambda_{t-}-z_n$.
Next, we let
\begin{eqnarray}
&& \sigma_1=\sigma_1(s,\Delta):=\sup\{\widetilde{s}\in[0,s]\!:\widetilde{s}\in\Delta\NN\}, \\
&& \sigma_2=\sigma_2(n,s,\Delta)
:=\inf\{\widetilde{s}\in[s,T]\!:\Lambda^\Delta_{\widetilde{s}}\le\Lambda_{t-}-z_n\}.
\end{eqnarray}
Then, by Definition \ref{def:Euler},
\begin{equation}\label{3terms}
\begin{split}
\frac{1}{d}\big((\Lambda^\Delta_{\sigma_1})^d-(\Lambda^\Delta_{\sigma_2})^d\big)
\leq&\,\int_{\RR_+\setminus [\Lambda^\Delta_{\sigma_2}-\varepsilon,\Lambda^\Delta_{\sigma_1}+\varepsilon]} u^\Delta(\sigma_2,x)\,\nu(\mathrm{d}x)
+\int_{(\Lambda^\Delta_{\sigma_2}-\varepsilon)^+}^{\Lambda^\Delta_{\sigma_1}+\varepsilon} u^\Delta(\sigma_2,x)\,\nu(\mathrm{d}x) \\
& -\int_{\RR_+} u^\Delta(\sigma_1,x)\,\nu(\mathrm{d}x).
\end{split}
\end{equation}



\smallskip

To study the first summand on the right-hand side of \eqref{3terms} we apply Definition \ref{def:Euler} and~get
\begin{equation*}
\begin{split}
u^\Delta(\sigma_2,x)
& =\E^x\big[\mathbf{1}_{\{\tau_{\Lambda^\Delta_{\sigma_2-\cdot}}<\sigma_2-\sigma_1\}}\,H(R_{\tau_{\Lambda^\Delta_{\sigma_2-\cdot}}})\big]
+\E^x\big[\mathbf{1}_{\{\tau_{\Lambda^\Delta_{\sigma_2-\cdot}}\ge \sigma_2-\sigma_1\}}\,u^\Delta(\sigma_1,R_{\sigma_2-\sigma_1})\big] \\
& \le H\big(\inf_{r\in[0,\sigma_2)} \Lambda^\Delta_r\big)
\,\PP^x(\tau_{\Lambda^\Delta_{\sigma_2-\cdot}}<\sigma_2-\sigma_1)
+\E^x[u^\Delta(\sigma_1,R_{\sigma_2-\sigma_1})],
\end{split}
\end{equation*}
for $x\in \RR_+\setminus [\Lambda^\Delta_{\sigma_2}-\varepsilon,\Lambda^\Delta_{\sigma_1}+\varepsilon]$.
~To analyze the second summand on the right-hand side of~\eqref{3terms} we note that 
\begin{equation*}
\begin{split}
u^\Delta(\sigma_2,x)
&= \E^x\big[\mathbf{1}_{\{\tau_{\Lambda^\Delta_{\sigma_2-\cdot}}<\sigma_2-\sigma_1\}}\,H(R_{\tau_{\Lambda^\Delta_{\sigma_2-\cdot}}})\big]
+\E^x\big[\mathbf{1}_{\{\tau_{\Lambda^\Delta_{\sigma_2-\cdot}}\ge \sigma_2-\sigma_1\}}\,u^\Delta(\sigma_1,R_{\sigma_2-\sigma_1})\big] \\
&\le H\big(x\wedge\inf_{r\in[0,\sigma_2)} \Lambda^\Delta_r\big)\vee 1
\end{split}
\end{equation*}
for all $x\in[(\Lambda^\Delta_{\sigma_2}-\varepsilon)^+,\Lambda^\Delta_{\sigma_1}+\varepsilon]$, whereas for $x\in[\Lambda^\Delta_{\sigma_2}+\varepsilon,\Lambda^\Delta_{\sigma_1}-\varepsilon]$,
\begin{equation*}
u^\Delta(\sigma_2,x)\le \big(H\big(x\wedge\inf_{r\in[0,\sigma_2)} \Lambda^\Delta_r\big)\vee 1\big)
\,\PP^x(\tau_{x-\varepsilon}\wedge\tau_{x+\varepsilon}<\sigma_2-\sigma_1)+H(x-\varepsilon)
\end{equation*}
(since for such $x$, it holds $\tau_{\Lambda^\Delta_{\sigma_2-\cdot}}<\sigma_2-\sigma_1$ on the event $\{\tau_{x-\varepsilon}\!\wedge\!\tau_{x+\varepsilon}\ge\sigma_2-\sigma_1\}$).
~Plugging the estimates of this paragraph into \eqref{3terms}, we arrive at
\begin{equation}\label{eq.JmpUpperBnd.mainProp.Lambda.sigma.1.2}
\begin{split}
&\;\frac{1}{d}\big((\Lambda^\Delta_{\sigma_1})^d-(\Lambda^\Delta_{\sigma_2})^d\big) \\
&\le \int_{\Lambda^\Delta_{\sigma_2}+\varepsilon}^{\Lambda^\Delta_{\sigma_1}-\varepsilon} H(x-\varepsilon)\,\nu(\mathrm{d}x)
+\int_{\RR_+\backslash[\Lambda^\Delta_{\sigma_2}-\varepsilon,\Lambda^\Delta_{\sigma_1}+\varepsilon]} \E^x[u^\Delta(\sigma_1,R_{\sigma_2-\sigma_1})]\,\nu(\mathrm{d}x)
-\int_{\RR_+} u^\Delta(\sigma_1,x)\,\nu(\mathrm{d}x) \\
&\;\;\;+ \bigg(4\varepsilon + \int_{\Lambda^\Delta_{\sigma_2}+\varepsilon}^{\Lambda^\Delta_{\sigma_1}-\varepsilon} \PP^x(\tau_{x-\varepsilon}\!\wedge\!\tau_{x+\varepsilon}<\sigma_2\!-\!\sigma_1)\,\mathrm{d}x\bigg)
\,\sup_{x\in[(\Lambda^\Delta_{\sigma_2}-\varepsilon)^+,\Lambda^\Delta_{\sigma_1}+\varepsilon]}x^{d-1}\big(H\big(x\!\wedge\!\inf_{r\in[0,\sigma_2)} \Lambda^\Delta_r\big)\!\vee\! 1\big)  \\
&\;\;\;+ \int_{\RR_+\backslash[\Lambda^\Delta_{\sigma_2}-\varepsilon,\Lambda^\Delta_{\sigma_1}+\varepsilon]}
H\big(\inf_{r\in[0,\sigma_2)} \Lambda^\Delta_r\big)
\,\PP^x(\tau_{\Lambda^\Delta_{\sigma_2-\cdot}}<\sigma_2-\sigma_1)\,\nu(\mathrm{d}x),
\end{split}
\end{equation}
for all sufficiently small $\Delta,\varepsilon>0$.


\medskip

Next, we take the limit $\Delta\downarrow0$ on both sides.~By construction, we have $\lim_{\Delta\downarrow0} \Lambda^\Delta_{\sigma_1}=\Lambda_s$. In addition, thanks to Lemma \ref{down_jumps_small.3}, we have $\underline{\lim}_{\Delta\downarrow0} \Lambda^\Delta_{\sigma_2} = \Lambda_{t-}-z_n-\delta(s)$, with $\delta(s)\geq0$ such that $\lim_{s\uparrow t} \delta(s)=0$. ~This, in particular, means that for all large enough $s<t$ and small enough $\varepsilon>0$ we have $\underline{\lim}_{\Delta\downarrow0} \Lambda^\Delta_{\sigma_2}-\varepsilon>0$.
~Using these observations, Proposition \ref{prop:tightness} and Remark \ref{rem:uDelta.bound} we take limits in \eqref{eq.JmpUpperBnd.mainProp.Lambda.sigma.1.2}, as $\Delta\downarrow0$, to obtain
\begin{equation*}
\begin{split}
&\;\frac{1}{d}\big((\Lambda_s)^d-(\Lambda_{t-}-z_n-\delta(s))^d\big)\\
&\le
\int_{\Lambda_{t-}-z_n-\delta(s)+\varepsilon}^{\Lambda_s-\varepsilon} H(x-\varepsilon)\,\nu(\mathrm{d}x)
-\int_{\RR_+} u(s,x)\,\nu(\mathrm{d}x)+C_1\varepsilon \\
&\;\;\;+\underset{\Delta\downarrow0}{\overline{\lim}}\, \int_{\RR_+\backslash[\Lambda^\Delta_{\sigma_2},\Lambda^\Delta_{\sigma_1}]} \E^x[u^\Delta(\sigma_1,R_{\sigma_2-\sigma_1})]\,\nu(\mathrm{d}x) 
\!+\! C_2\,\underset{\Delta\downarrow0}{\overline{\lim}}\,
\int_{\Lambda^\Delta_{\sigma_2}+\varepsilon}^{\Lambda^\Delta_{\sigma_1}-\varepsilon} \PP^x(\tau_{x-\varepsilon}\!\wedge\!\tau_{x+\varepsilon}\!<\!\sigma_2\!-\!\sigma_1)\,\nu(\mathrm{d}x) \\
&\;\;\;+C_3\,\underset{\Delta\downarrow0}{\overline{\lim}}\, \int_{\RR_+\backslash[\Lambda^\Delta_{\sigma_2}-\varepsilon,\Lambda^\Delta_{\sigma_1}+\varepsilon]} \PP^x(\tau_{\Lambda^\Delta_{\sigma_2-\cdot}}<\sigma_2-\sigma_1)\,\nu(\mathrm{d}x).
\end{split}
\end{equation*}

\smallskip

Finally, we pass to the limit in the above, as $s\uparrow t$, to deduce 
\begin{equation}\label{eq.JmpUpperBnd.beforeStep2}
	\begin{split}
		& \frac{1}{d}\big((\Lambda_{t-})^d-(\Lambda_{t-}-z_n)^d\big) 
		\leq
		\int_{\Lambda_{t-}-z_n+\varepsilon}^{\Lambda_{t-}-\varepsilon} H(x-\varepsilon)\,\nu(\mathrm{d}x)
	    -\int_{\Lambda_{t-}-z_n}^{\Lambda_{t-}} u(t-,x)\,\nu(\mathrm{d}x)	
	    +C_1\varepsilon \\
		& \qquad\qquad\quad\; -\int_{\RR_+\backslash[\Lambda_{t-}-z_n,\Lambda_{t-}]} u(t-,x)\,\nu(\mathrm{d}x)	
		+\underset{s\uparrow t}{\overline{\lim}}\;\underset{\Delta\downarrow0}{\overline{\lim}}\, \int_{\RR_+\backslash[\Lambda^\Delta_{\sigma_2},\Lambda^\Delta_{\sigma_1}]} \E^x[u^\Delta(\sigma_1,R_{\sigma_2-\sigma_1})]\,\nu(\mathrm{d}x) \\
		& \qquad\qquad\quad\; +C_2\,\underset{s\uparrow t}{\overline{\lim}}\;\underset{\Delta\downarrow0}{\overline{\lim}}\,
		\int_{\Lambda^\Delta_{\sigma_2}+\varepsilon}^{\Lambda^\Delta_{\sigma_1}-\varepsilon} \PP^x(\tau_{x-\varepsilon}\!\wedge\!\tau_{x+\varepsilon}<\sigma_2-\sigma_1)\,\nu(\mathrm{d}x) \\
		& \qquad\qquad\quad\; +C_3\,\underset{s\uparrow t}{\overline{\lim}}\;\underset{\Delta\downarrow0}{\overline{\lim}}\, \int_{\RR_+\backslash[\Lambda^\Delta_{\sigma_2}-\varepsilon,\Lambda^\Delta_{\sigma_1}+\varepsilon]} \PP^x(\tau_{\Lambda^\Delta_{\sigma_2-\cdot}}\!<\!\sigma_2-\sigma_1)\,\nu(\mathrm{d}x).
	\end{split}
\end{equation}
At this point, it suffices to show that the second, the third, and the fourth lines in the latter display are non-positive.~Indeed, once the latter is established, we can let $\varepsilon\downarrow0$ and recall the definition of $z_n$, to obtain
\begin{align*}
\frac{1}{d}\big((\Lambda_{t-})^d-(\Lambda_{t-}-z_n)^d\big) 
&\leq \int_{\Lambda_{t-}-z_n}^{\Lambda_{t-}} H(x)\,\nu(\mathrm{d}x)
-\int_{\Lambda_{t-}-z_n}^{\Lambda_{t-}} u(t-,x)\,\nu(\mathrm{d}x) \\
&\leq \frac{1}{d}\big((\Lambda_{t-})^d-(\Lambda_{t-}-z_n)^d\big) - \frac{1}{n},\\
\end{align*}
which is the desired contradiction.

\medskip

The fourth line in \eqref{eq.JmpUpperBnd.beforeStep2},
\begin{equation*}
\begin{split}
& \;C_3\,\underset{s\uparrow t}{\overline{\lim}}\;\underset{\Delta\downarrow0}{\overline{\lim}}\, \int_{\RR_+\backslash[\Lambda^\Delta_{\sigma_2}-\varepsilon,\Lambda^\Delta_{\sigma_1}+\varepsilon]} \PP^x(\tau_{\Lambda^\Delta_{\sigma_2-\cdot}}\!<\!\sigma_2-\sigma_1)\,\nu(\mathrm{d}x) \\
& = C_3\,\underset{s\uparrow t}{\overline{\lim}}\;\underset{\Delta\downarrow0}{\overline{\lim}}\, \int_{\RR_+\backslash[\Lambda_{t-}-z_n-\varepsilon,\Lambda_s+\varepsilon]} \PP^x(\tau_{\Lambda^\Delta_{\sigma_2-\cdot}}\!<\!\sigma_2-\sigma_1)\,\nu(\mathrm{d}x),
\end{split}
\end{equation*}
vanishes because $\Lambda^\Delta_{\sigma_2-\cdot}|_{[0,\sigma_2-\sigma_1]}$ takes values in $[\Lambda_{t-}-z_n-\varepsilon/2,\Lambda_s+\varepsilon/2]$ for all $s<t$ close enough to $t$ and all $\Delta>0$ small enough; $\lim_{\Delta\downarrow0}\,(\sigma_2-\sigma_1)=t-s$; the process $R^x$ dominates a standard Brownian motion started from $x$ for $x>\Lambda_s+\varepsilon$; and the process $(R^x_{\cdot\wedge\tau_{\Lambda_{t-}-z_n-\varepsilon/2}})^2$ can be compared to a Brownian motion with drift $d$ and diffusion coefficient $2(\Lambda_{t-}-z_n-\varepsilon/2)$ started from $x^2$ for $x<\Lambda_{t-}-z_n-\varepsilon$ (see the second paragraph in the proof of Lemma \ref{down_jumps_small} for a very similar argument). 

\medskip

To see that the third line in \eqref{eq.JmpUpperBnd.beforeStep2},
\begin{equation*}	
\begin{split}
&\; C_2\,\underset{s\uparrow t}{\overline{\lim}}\;\underset{\Delta\downarrow0}{\overline{\lim}}\,
\int_{\Lambda^\Delta_{\sigma_2}+\varepsilon}^{\Lambda^\Delta_{\sigma_1}-\varepsilon} \PP^x(\tau_{x-\varepsilon}\!\wedge\!\tau_{x+\varepsilon}<\sigma_2-\sigma_1)\,\nu(\mathrm{d}x) \\
& \le C_2\,\underset{s\uparrow t}{\overline{\lim}}\;\underset{\Delta\downarrow0}{\overline{\lim}}\,
\int_{\Lambda_{t-}-z_n+\varepsilon}^{\Lambda_s-\varepsilon} \PP^x(\tau_{x-\varepsilon}<\sigma_2-\sigma_1)+\PP^x(\tau_{x+\varepsilon}<\sigma_2-\sigma_1)
\,\nu(\mathrm{d}x),
\end{split}
\end{equation*}
equals to $0$, we recall that $\lim_{\Delta\downarrow0}\,(\sigma_2-\sigma_1)=t-s$, dominate $R^x$ by a standard Brownian motion started from $x$ to bound $\PP^x(\tau_{x-\varepsilon}<\sigma_2-\sigma_1)$, and compare $(R^x_{\cdot\wedge\tau_{x+\varepsilon}})^2$ to a Brownian motion with drift $d$ and diffusion coefficient $2\Lambda_s$ started from $x^2$ to estimate $\PP^x(\tau_{x+\varepsilon}\!<\!\sigma_2\!-\!\sigma_1)$.

\medskip

Lastly, we consider the term
\begin{equation}\label{last_term}
\underset{s\uparrow t}{\overline{\lim}}\;\underset{\Delta\downarrow0}{\overline{\lim}}\, \int_{\RR_+\backslash[\Lambda^\Delta_{\sigma_2},\Lambda^\Delta_{\sigma_1}]} \E^x[u^\Delta(\sigma_1,R_{\sigma_2-\sigma_1})]\,\nu(\mathrm{d}x)
-\int_{\RR_+\backslash[\Lambda_{t-}-z_n,\Lambda_{t-}]} u(t-,x)\,\nu(\mathrm{d}x).
\end{equation}
Applying Lemma \ref{down_jumps_small.3} and the Dominated Convergence Theorem (recall Remark \ref{rem:uDelta.bound}) we see
\begin{equation*}
\begin{split}
&\;\,\underset{\Delta\downarrow0}{\overline{\lim}}\, \bigg|\int_{\RR_+\backslash[\Lambda^\Delta_{\sigma_2},\Lambda^\Delta_{\sigma_1}]} \E^x[u^\Delta(\sigma_1,R_{\sigma_2-\sigma_1})]\,\nu(\mathrm{d}x) 
-\int_{\RR_+\backslash[\Lambda_{t-}-z_n,\Lambda_s]} \E^x[u^\Delta(\sigma_1,R_{\sigma_2-\sigma_1})]\,\nu(\mathrm{d}x)
\bigg| \\
&\leq \underset{\Delta\downarrow0}{\overline{\lim}}\,\int_{[\Lambda_{t-}-z_n-\delta(s),\Lambda_{t-}-z_n]} \E^x[u^\Delta(\sigma_1,R_{\sigma_2-\sigma_1})]\,\nu(\mathrm{d}x)
\leq C_4\,\delta(s).
\end{split}
\end{equation*}
Using the Dominated Convergence Theorem and the definition of the Bessel process we get
\begin{align*}
&\;\int_{\RR_+\backslash[\Lambda_{t-}-z_n,\Lambda_s]} \E^x[u^\Delta(\sigma_1,R_{\sigma_2-\sigma_1})]\,\nu(\mathrm{d}x)\\
&= \frac{\Gamma(d/2+1)}{d\pi^{d/2}}\,\E\bigg[\int_{\RR^d} \bone_{\RR_+\backslash[\Lambda_{t-}-z_n,\Lambda_s]}(|x|)
\, u^\Delta(\sigma_1,|x+B_{\sigma_2-\sigma_1}|)\,\mathrm{d}x\bigg]\\
&= \frac{\Gamma(d/2+1)}{d\pi^{d/2}}\,\E\bigg[\int_{\RR^d} \bone_{\RR_+\backslash[\Lambda_{t-}-z_n,\Lambda_s]}(|x-B_{\sigma_2-\sigma_1}|)
\, u^\Delta(\sigma_1,|x|)\,\mathrm{d}x\bigg],
\end{align*}
where $\Gamma$ is the Gamma function and $B$ is a standard Brownian motion in $\RR^d$.~Then, Lemma~\ref{lem:u_conv} and the Dominated Convergence Theorem yield
\begin{align*}
&\;\underset{\Delta\downarrow0}{\overline{\lim}} \int_{\RR_+\backslash[\Lambda_{t-}-z_n,\Lambda_s]} \E^x[u^\Delta(\sigma_1,R_{\sigma_2-\sigma_1})]\,\nu(\mathrm{d}x)\\
&= \frac{\Gamma(d/2+1)}{d\pi^{d/2}}\,\E\bigg[\int_{\RR^d} \bone_{\RR_+\backslash[\Lambda_{t-}-z_n,\Lambda_s]}(|x-B_{t-s}|)
\, u(s,|x|)\,\mathrm{d}x\bigg].
\end{align*}
Taking $s\uparrow t$ we obtain
\begin{align*}
&\;\underset{s\uparrow t}{\overline{\lim}}\;\underset{\Delta\downarrow0}{\overline{\lim}}\, \int_{\RR_+\backslash[\Lambda_{t-}-z_n,\Lambda_s]} \E^x[u^\Delta(\sigma_1,R_{\sigma_2-\sigma_1})]\,\nu(\mathrm{d}x)\\
&= \frac{\Gamma(d/2+1)}{d\pi^{d/2}} \int_{\RR^d} \bone_{\RR_+\backslash[\Lambda_{t-}-z_n,\Lambda_{t-}]}(|x|)
\, u(t-,|x|)\,\mathrm{d}x
 = \int_{\RR_+\backslash[\Lambda_{t-}-z_n,\Lambda_{t-}]} u(t-,x)\,\nu(\mathrm{d}x).
\end{align*}
Collecting the above four displays we conclude that the expression in \eqref{last_term} equals to $0$.
\qed

\bigskip\bigskip

\bibliographystyle{amsalpha}
\bibliography{Main}

\bigskip\bigskip\medskip

\end{document}